\documentclass[notitlepage, 11pt]{article}

\usepackage{color}
\usepackage{latexsym}
\usepackage{amsmath, amssymb, amsfonts}

\DeclareMathOperator*{\argmin}{arg\,min}

\newtheorem{lemma}{Lemma}
\newtheorem{theorem}{Theorem}

\newtheorem{definition}{Definition}

\newtheorem{proposition}{Proposition}

\def\thelemma{\arabic{section}.\arabic{lemma}}
\def\thetheorem{\arabic{section}.\arabic{theorem}}
\def\thecorollary{\arabic{section}.\arabic{corollary}}
\def\thedefinition{\arabic{section}.\arabic{definition}}
\def\theexample{\arabic{section}.\arabic{example}}
\def\theproposition{\arabic{section}.\arabic{proposition}}

\def\theassumption{\arabic{section}.\arabic{assumption}}

\def\theremark{\arabic{section}.\arabic{remark}}

\newcommand{\manualnames}[1]{
\def\thelemma{#1.\arabic{lemma}}
\def\thetheorem{#1.\arabic{theorem}}
\def\thecorollary{#1.\arabic{corollary}}
\def\thedefinition{#1.\arabic{definition}}
\def\theexample{#1.\arabic{example}}
\def\theproposition{#1.\arabic{proposition}}
\def\theassumption{#1.\arabic{assumption}}
\def\theremark{#1.\arabic{remark}}
}

\newcommand{\beginsec}{
\setcounter{lemma}{0}
\setcounter{theorem}{0}
\setcounter{corollary}{0}
\setcounter{definition}{0}
\setcounter{example}{0}
\setcounter{proposition}{0}
\setcounter{condition}{0}
\setcounter{assumption}{0}
\setcounter{conjecture}{0}
\setcounter{problem}{0}
\setcounter{remark}{0}
}

\newcommand{\noi}{\noindent}

\newcommand{\E}{\mathbb{E}}
\newcommand{\R}{\mathbb{R}}

\newcommand{\N}{\mathbb{N}}
\newcommand{\I}{\mathcal{I}}
\newcommand{\Ir}{\mathbb{I}}

\newcommand{\la}{\lambda}
\newcommand{\sig}{\sigma}

\newcommand{\eps}{\varepsilon}

\newcommand{\Barr}{{\beta_0}}
\newcommand{\taud}{\tau_\delta}

\newcommand{\tauald}{\tau_{\bar\al,\delta}}
\newcommand{\phiu}{\phi_u}
\newcommand{\phiv}{\phi_v}
\newcommand{\cc}{c}
\newcommand{\cco}{c_1}
\newcommand{\cct}{c_2}
\newcommand{\ye}{y_\eps}
\newcommand{\xe}{x_\eps}
\newcommand{\xd}{x_\delta}
\newcommand{\gd}{g_\delta}
\newcommand{\ld}{l_\delta}

\newcommand{\ph}{\varphi}
\newcommand{\vr}{\varrho}
\newcommand{\al}{\alpha}

\newcommand{\del}{\delta}
\newcommand{\om}{\omega}

\newcommand{\Gam}{\mathnormal{\Gamma}}

\newcommand{\calC}{{\cal C}}
\newcommand{\calD}{{\cal D}}
\newcommand{\calE}{{\cal E}}

\newcommand{\calH}{{\cal H}}
\newcommand{\calI}{{\cal I}}

\newcommand{\calL}{{\cal L}}

\newcommand{\calP}{{\cal P}}
\newcommand{\calQ}{{\cal Q}}

\newcommand{\AC}{{\cal AC}}

\newcommand{\oo}{\overline}
\newcommand{\skp}{\vspace{\baselineskip}}

\newcommand{\w}{\wedge}

\newcommand{\iy}{\infty}

\newcommand{\ds}{\displaystyle}
\newcommand{\A}{{\cal A}}

\newcommand{\IA}{{\it IA}}
\newcommand{\ST}{{\it ST}}
\newcommand{\qed}{\hfill $\Box$}

\newcommand{\uh}{\underline{h}}
\newcommand{\ur}{\underline{r}}
\newcommand{\uq}{\underline{q}}
\newcommand{\uy}{\underline{y}}
\newcommand{\ux}{\underline{x}}
\newcommand{\uc}{\underline{c}}

\newcommand{\upsi}{\underline{\psi}}
\newcommand{\uxi}{\underline{\xi}}
\newcommand{\utheta}{\underline{\theta}}
\newcommand{\uzeta}{\underline{\zeta}}
\newcommand{\ual}{\underline{\alpha}}

\newcommand{\uph}{\underline{\varphi}}
\newcommand{\uvr}{\underline{\varrho}}

\newcommand{\Jr}{\mathbb{J}}

\title{A differential game for a multiclass queueing model
in the moderate-deviation heavy-traffic regime\thanks{Research
supported in part by the ISF (Grant 1315/12)
and the Technion fund for promotion of research}}
\author{Rami Atar\thanks{
Department of Electrical Engineering,
Technion--Israel Institute of Technology,
Haifa 32000, Israel}
\and
Asaf Cohen\thanks{Department of Mathematics,
University of Michigan,
Ann Arbor, MI 48109, USA}
}

\date{October 14, 2014}

\begin{document}

\maketitle

\begin{abstract}
We study a differential game that governs the moderate-deviation
heavy-traffic asymptotics
of a multiclass single-server queueing control problem
with a risk-sensitive cost.
We consider a cost set on a finite
but sufficiently large time horizon, and show
that this formulation leads to stationary feedback policies for the game.
Several aspects of the game are explored, including its characterization
via a (one-dimensional) free boundary problem,
the semi-explicit solution of an optimal strategy, and the specification
of a saddle point. We emphasize the analogy to the well-known
Harrison-Taksar free boundary problem which plays a similar role
in the diffusion-scale heavy-traffic literature.

\skp

\noi{\bf AMS subject classifications:}\,\, 49N70, 60F10, 60K25, 93E20

\skp

\noi{\bf Keywords:}\,\,
multi-class single-server queue, moderate deviations,
heavy traffic, risk-sensitive cost,
the Harrison-Taksar free boundary problem
\end{abstract}

\section{Introduction}
\beginsec

This paper is concerned with the moderate-deviation (MD)
scale heavy-traffic analysis of a queueing
control problem. The model that is treated consists of a server that devotes its time to
customers from a number of classes, where a decision maker (DM)
controls the allocation of the server's effort to each of the classes.
Customers of each class are kept in a buffer of
finite length and so new arrivals of customers of a given class
are lost when the corresponding buffer is full.
The DM may reject arrivals even when buffers are not full.
The cost, that is of risk-sensitive type and is rescaled at the MD
regime, accounts for holding of customers in the buffers as well as for rejections.
The term `heavy traffic' refers to the imposition of
a critical load condition, of the traffic intensity being close to one.

Heavy traffic analysis is traditionally carried out under the regime
of {\it diffusion-scale deviations}
(sometimes also referred to as {\it ordinary deviations}),
where a vast variety of queueing control problems have been considered.
However, there are very few results on the corresponding MD regime.
This paper aims at continuing the direction started in \cite{ata-bis}
to develop an approach to queueing control problems at the MD scale.
While \cite{ata-bis} finds various similarities to the diffusion limit theory,
a major dissimilarity between the two formulations
is at the heart of the present paper, as described below.

In \cite{ata-dup-shw} it has been shown that the
{\it large deviation} asymptotics of risk-sensitive control problems associated with
queueing models are governed by certain zero-sum differential games (DG)
(this approach was inspired by and built on analogous
treatments of small noise controlled diffusion models, such as \cite{fle-sou}).
In these games, the
dynamics are given by a controlled ordinary differential equation (ODE),
and the players are the DM and nature, that attempt to minimize and, respectively,
maximize a cost. The DM
again controls server's effort allocation and admission/rejection (except that
the model no longer accounts for customers and so terms such as
effort and admission/rejection are interpreted
as suitable elements of the ODE), and nature may perturb the arrival and
service rates. We will sometimes refer to the DM and nature as the minimizing
player and maximizing player, respectively, or simply as minimizer and maximizer.
The game's cost consists of the original holding and rejection penalties,
and an additional term associated with the perturbation of the aforementioned rates.
The latter term is determined by the underlying large deviation action functional.
The treatment of a queueing control problem in the MD scale in
\cite{ata-bis} is similar to that of large deviations as far as the tools are concerned,
but it is argued in that paper that the games one obtains in the MD
regime tend to be solvable explicitly or semi-explicitly more often than
under large deviation scaling.
In fact, the specific DG identified in \cite{ata-bis} turned out to be
solvable due to a property of the model
referred to in the heavy traffic literature as {\it pathwise minimality},
that enables approaching the game in a straightforward manner
without appealing to dynamic programming methods.
However, this property is quite special and it is therefore desired to find methods
that apply in greater generality.
In the present paper we aim at models that do not possess
the pathwise minimality property, for which the solution methods must rely on
tools such as dynamic programming.
The model presented above is a prototype for
models where a dynamic programming approach is required; specifically,
the Bellman equation for the model
identifies a {\it free boundary point} that is a crucial ingredient
of the game's solution.
The paper is devoted to the derivation and study of the DG itself, whereas
the convergence of the risk-sensitive value to that of the DG is the subject of
a work in progress \cite{atar-cohen}.

As mentioned above, it was argued in \cite{ata-bis} that unlike treatment
at the large deviations regime, MD analysis
parallels that of the diffusion regime in several ways.
The present paper's contribution is best explained in the context of this relation.
While the asymptotics of the MD scaled
risk-sensitive control problem leads to a DG,
for diffusion scaled (risk-neutral) control problems, they
are typically described in terms of a {\it Brownian control problem} (BCP),
namely a control problem involving diffusion processes.
The BCP and the DG are similar in their structure,
where the difference lies in that in the BCP, uncertainty is captured
by the stochasticity of the Brownian motion term,
whereas in the DG the role of uncertainty is played by an adversarial disturbance,
in the form of the arrival and service rates being controlled by nature.
The structural similarity is manifested by the fact
that the mapping from the Brownian motion
paths to the optimal control for the BCP is exactly the one used for
the DG's optimal strategy.
Moreover, a phenomenon known as {\it state space collapse}, where a BCP
formulated in dimension $d$ reduces to a 1-dimensional control problem,
occurs in the corresponding DG as well (see \cite{ata-bis} for details).

In the situation of the present paper, an analogy
is less obvious, for reasons having to do with a disadvantage of
working with a discount factor in a risk-sensitive cost formulation.
In the case of risk-neutral cost, a standard Markov decision process formulation based on
discounted cost has the advantage over
a finite time horizon cost, that optimal controls are
given by stationary feedback. (For a typical queueing model this would often mean
that decisions such as admission, server effort allocation and routing
depend on the collection of queue lengths of the system in way that
does not change over time.) However, it is well understood
that a discounted version of a risk-sensitive
cost leads to optimal feedback that is {\it time dependent}
(this has been noticed first in \cite{chung}; several aspects of this
issue were treated in \cite{marcus}). Thus, although it
is desirable to aim at a formulation that leads to stationary feedback
because of the elegant structure of the solution,
there does not seem to be a way to do so based on
a discounted risk-sensitive cost. We address this issue by considering
a risk-sensitive cost over a finite, but sufficiently
large time horizon, in a formulation
that does lead to stationary optimal feedback controls as far as the DG is concerned.

It is instructive to notice that thanks to the desired property alluded to above,
the aforementioned analogy goes a step further.
The {\it Harrison-Taksar free boundary problem} addresses a 1-dimensional
BCP with an additive singular control term, for which
the solution is given as a reflecting Brownian motion on an interval $[a,b]$,
where $a$ and $b$ form a free boundary in a suitable Bellman
equation \cite{har-tak}. The significance of the free boundary point $b$ in applications to
diffusion scaled queueing models is that it determines a threshold, in queue length
terms, for rejection of arrivals.
A (distinct) free boundary point identified in the present paper, denoted by
$\beta_0$, plays the same role in the
MD scale treatment: jobs are to be rejected once this point is exceeded.
As far as the DG is concerned,
an optimal strategy for the DM is one that applies the Skorohod map
on the interval $[0,\beta_0]$ to the dynamics of the game, regardless of the
behavior of nature. Note the similar role played by the Skorohod map in both
scenarios.

As far as analyzing the DG is concerned, the main contribution is the identification
of the optimal strategy alluded to above. Apart from that,
we state and prove several additional results, some of which form the
basis for the result regarding the optimal strategy. Specifically, we find an equivalent
game that is easier to analyze, characterize the DG's value function
as the unique viscosity solution to the underlying Bellman equation, solve the equation
explicitly, and identify a saddle point for the game as well as the optimal behavior
of nature under the optimal play of the DM.

To summarize the contribution of this paper, we have
\begin{itemize} \itemsep0em
\item
 identified a queueing control problem
 that is prototypical to situations where a dynamic programming
 approach is required to treat the DG governing the MD heavy traffic limit,
\item proposed a risk-sensitive cost
 that leads to optimal stationary feedback for the DG,
\item solved the Bellman equation, found an optimal strategy for the DM
 and studied properties of the DG.
\end{itemize}

Let us finally mention that the model studied here was recently
treated in diffusion-scale heavy-traffic \cite{ata-shi}.

The organization of the paper is as follows. In Section \ref{sec2} we
describe the queueing model, introduce the MD scaling, present our approach
to formulation of the cost and introduce the DG.
In Section \ref{sec_3} we analyze the DG by introducing the hitting time game,
present the Bellman equation, characterize the
DG's value function as its unique viscosity solution,
present an explicit optimal strategy for the DM,
as well as a saddle point result.

We use the following notation.
For a positive integer $k$ and $a,b\in\R^k$, $a\cdot b$ denotes the usual
scalar product. We denote $[0, \iy)$ by $\R_+$.
Denote by $\AC([0,T],\R^k)$, $\calC([0,T],\R^k)$ and $\calD([0,T], \R^k)$ the
spaces of absolutely continuous functions [resp.,
continuous functions, functions that are
right-continuous with finite left limits (RCLL)] mapping $[0,T]\to\R^k$.
Write $\AC_0([0,T],\R^k)$ and $C_0([0,T],\R^k)$ for the subsets of the corresponding
function spaces, of functions that start at zero.
Endow the space $\calD([0,T],\R^k)$ with the usual Skorohod topology.

\section{Model and results}\label{sec2}
\beginsec

\subsection{The single-class model}\label{sec2a}

We describe the probabilistic model that provides the main motivation
for our DG formulation. The DG is derived from this model by considering
scaling limits of the latter. Because the focus of the present paper is on the analysis
of the DG itself, the derivation will only be carried out in
a heuristic manner. In a work in progress \cite{atar-cohen}
this relation is explored further and a rigorous relation is established.

We first present a basic model, then an extension of it.
In the basic model, that we will call the {\it single-class model},
customers arrive into a buffer with finite room and are served by a single server
in the order of arrival. The decision maker may reject arrivals even if the buffer
is not full, but in order to keep the buffer length constraint it must
always reject arrivals that occur when the buffer is full.
The cost functional will involve penalties based on queue-length as well as rejection.

Scaling is introduced by considering a sequence of systems, indexed by $n\in\N$.
In the $n$-th system, arrivals occur according to a renewal process, $A^n$, for which the mean
inter-arrival time is given by $1/\la^n$.
Denote by $S^n(t)$ the
number of service completions by the time the server has worked for $t$ units of time.
Then the inter-jump times of this process correspond to service time durations,
and it is assumed that $S^n$ is a renewal process (that is, the service times
are i.i.d.). It is also assumed that $A^n$ and $S^n$ are mutually independent.
The mean service time is denoted by $1/\mu^n$.
It is assumed that the squared coefficient of variation for
inter-arrival and service time distribution are
fixed, and these constants are
denoted by $\sig^2_{\IA}\in (0,\infty)$ and $\sig^2_{\ST}\in (0,\infty)$,
respectively.

Denoting the number of customers in the $n$-th system at time $t$ by $X^n(t)$,
the number of rejections until time $t$ by $R^n(t)$, and
the cumulative business time till $t$ by $T^n(t)$, we have the balance equation
\begin{equation}\label{eq008}
X^n(t)=X^n(0)+A^n(t)-S^n(T^n(t))-R^n(t),\qquad t\ge 0.
\end{equation}
The process $U^n=(T^n, R^n)$ is regarded as a control.
The finite buffer constraint reads $X^n(t)\le D^n$ for all $t$, where $D^n$
are given constants.

We consider the {\it moderate deviations rate parameters} $\{b_n\}$, that form a sequence,
fixed throughout, with the property that $\lim b_n=\infty$ while $\lim n^{-1/2}b_n=0$, as
$n\to\iy$.
We assume that as $n\to\iy$,
\[
\la^n=\la n+\tilde\la b_n\sqrt{n}+o(b_n\sqrt{n}),
\qquad
\mu^n=\mu n+\tilde\mu b_n\sqrt{n}+o(b_n\sqrt{n}),
\]
where $\la\in(0,\iy)$, $\mu\in(0,\iy)$, $\tilde\la\in\R$ and $\tilde\mu\in\R$ are fixed.
Moreover, it is assumed that the system is critically loaded, namely $\la =\mu$.

For $T>0$, let $\Ir_i(T,\cdot)$, $i=1,2$, be functions mapping
$\calD([0,T],\R)$ to $\R_+\cup\{+\iy\}$ defined as
\begin{equation}\label{eq015}
\Ir_i(T,\psi)=\begin{cases}
                \ds
                c_i\int_0^T\dot\psi^2(s)ds &\mbox{if}\ \psi\in\AC_0([0,T],\R),
\\ \\
                +\infty & \mbox{otherwise},
              \end{cases}
\end{equation}
where
$$
\cco = \frac{1}{2\mu\sig^2_{IA}} \text{ and } \cct=\frac{1}{2\mu\sig^2_{ST}}.
$$
In the sequel we will also use the notation
$$
\cc=(c_1^{-1}+c_2^{-1})^{-1}=\frac{1}{2\mu(\sig^2_{ST} + \sig^2_{AI})}.
$$
Let $\Ir(T,\psi):=\Ir_1(T,\psi^1)+\Ir_2(T,\psi^2)$.
Under suitable exponential moment conditions, the centered, scaled versions,
\[
  \tilde{A}^n(t)=\frac{1}{b_n\sqrt{n}}(A^n(t)-\la^n t),
  \quad \tilde{S}^n(t)=\frac{1}{b_n\sqrt{n}}(S^n(t)-\mu^n t),\quad t\ge0,
\]
of the renewal processes alluded to above satisfy a {\it moderate deviation principle}
\cite{Puhal-1999}, \cite{Puhal-Whitt}.
That is,
given $T>0$ and a bounded, continuous function $H:\calD([0,T],\R)^2\to\R$,
one has
\begin{equation}\label{03}
\lim_{n\to\iy}\frac{1}{b_n^2}\log E e^{b_n^2H(\tilde A^n,\tilde S^n)}
=\sup \{H(\psi)-\Ir(\psi): \psi\in\AC_0([0,T],\R)^2\}.
\end{equation}

Introduce scaled versions of the remaining two processes, namely
\begin{align}\label{eq010}
  \tilde{R}^n(t):=\frac{1}{b_n\sqrt{n}}R^n(t),
  \quad \tilde{X}^n(t)&:=\frac{1}{b_n\sqrt{n}}X^n(t).
\end{align}
The buffer size is assumed to be given by $D^n=b_n\sqrt n D$ for a constant $D>0$,
and so the buffer constraint can now be written as
\begin{equation}
\label{eq013}
\tilde X^n(t)\in[0,D],\qquad t\ge0.
\end{equation}

To define the cost, let $C^n(t)=\int_0^th(\tilde X^n(s))ds+r\tilde R^n(t)$,
where, throughout, $h$ is a continuous increasing function
from $[0,D]$ to $\R_+$ that satisfies $h(0) = 0$
and accounts for holding cost, and $r>0$ is a constant representing per-customer
rejection penalty. One may consider several versions of risk-sensitive costs based on
the process $C^n$, the most natural ones being
a cost defined over a finite time horizon,
\[
J_1=\frac{1}{b_n^2}\log E e^{b_n^2\int_t^T C^n(s)ds},
\]
a long-run average cost,
\[
J_2=\limsup_{T\to\iy}\frac{1}{T}\frac{1}{b_n^2}\log E e^{b_n^2\int_0^T C^n(t)dt},
\]
and two versions based on discount ($\kappa>0$ being the discount rate)
\[
J_3=\frac{1}{b_n^2}\log E e^{b_n^2\int_0^\iy e^{-\kappa t} C^n(t)dt},\qquad
J_4=\frac{1}{b_n^2}\log E\int_0^\iy e^{-\kappa t} e^{b_n^2 C^n(t)}dt.
\]
Because our setting is not Markovian, dynamic programming methods are not applicable
in a direct fashion for the $n$th system. They are only applicable as far as the
limit behavior is concerned. Thus the discussion that follows regards only
the limiting behavior as well as the Markovian setting, that occurs in the
special case when the driving renewal processes are Poisson.
Working with a cost defined over a finite time horizon, such as $J_1$,
suffers from the drawback
that it leads to optimal feedback policies that are non-stationary.
Long-run average cost of risk-sensitive type has been studied in several
papers starting from \cite{fle-mce}.
It leads to stationary optimal feedback
policies. However, working with this cost under scaling limits might
be quite demanding technically, and therefore we will not focus on this cost here.
Regarding the version $J_3$ of a discounted cost, it is known to lead to non-stationary
optimal feedback \cite{chung}. The reader is referred to \cite{marcus} for discussion
and several contributions regarding this issue.
The formulation we propose is close in spirit to $J_4$. Specifically,
given $n$ and $T$, we consider the following
cost associated with a control $U^n$, namely
\begin{align}\label{eq019}
\notag
J^n(U^n,T)&:=\frac{1}{b_n^2}
\log\E\Big[\int_0^T e^{b^2_nC^n(t)}dt\Big]
\\
&=\frac{1}{b_n^2}
\log\E\Big[\int_0^T e^{b^2_n[\int_0^t h(\tilde{X}^n(s))ds+r\tilde{R}^n(t))]}dt\Big],
\end{align}
where $(\tilde X^n,\tilde R^n)$ are the scaled versions
of the processes $(X^n,R^n)$ corresponding to the control $U^n$.
The value is given by
\begin{equation}\label{04}
V^n(T)=\inf_{U^n} J^n(U^n,T),
\end{equation}
and our aim will be to understand the limit
\begin{equation}
  \label{04-}
  \limsup_{T\to\iy}\limsup_{n\to\iy}V^n(T).
\end{equation}

\subsection{Differential game setting}\label{sec2b}

We provide a heuristic argument for the claim that the limit of $V^n$ is given
by a game, before turning to a precise definition of the game.
Note first that if one could write $C^n$ as function $H$ of the data
$(\tilde A^n,\tilde S^n)$ and the control $U^n$, namely $C^n=H(U^n,\tilde A^n,\tilde S^n)$,
then based on \eqref{03}, and ignoring questions of boundedness and continuity,
one would obtain in the limit a game, defined in terms of the cost
$H(u,\psi)-\Ir(\psi)$, maximized over $\psi$ (as in \eqref{03})
and minimized over $u$ (as in \eqref{04}).
Although this conclusion is not solid due to questions of continuity,
one can argue by appealing to specific
aspects of the dependence of $C^n$ on the data and control.
Specifically, one notes from \eqref{eq008} that
\begin{equation}
 \tilde{X}^n(t)=\tilde{X}^n(0) + \tilde y^nt
 +\tilde{A}^n(t)-\tilde{S}^n(T^n(t))
 + \tilde{Z}^n(t) - \tilde{R}^n(t),
 \label{eq011}
\end{equation}
where we denote
\begin{equation}\label{eq012}
\tilde{Z}^n(t):=\frac{\mu^n}{n}\frac{\sqrt{n}}{b_n}(t-T^n(t)),
\qquad \tilde y^n:=\tilde\la^n-\tilde\mu^n.
\end{equation}
Moreover, the processes $\tilde{Z}^n$ and $\tilde{R}^n$ are nondecreasing, whereas,
by assumption, $\tilde X^n(0)\to x$ and $\tilde y^n\to y$.
It is thus reasonable to consider a DG
with dynamics $\ph(t)=x+yt+\psi^1(t)-\psi^2(t)
+\zeta(t)-\varrho(t)$ where $\zeta$ and $\varrho$ are nondecreasing (corresponding
to the processes $Z^n$ and $R^n$), $\ph$ takes
values in $[0,D]$ (because of \eqref{eq013}), and the cost takes
the form
\[
\sup_{T\in[0,\iy)}\Big\{\int_0^Th(\ph(t))dt+
r\varrho(T)-\Ir(T,\psi)\Big\}.
\]
It is maximized over $\psi$ and minimized over
$u=(\varrho,\zeta)$. The precise formulation along these lines is given next.

Let $y=\tilde\la-\tilde\mu$.
Denote by $\calP=\calC_0([0,\iy),\R)$ and
\[
\calE=\{\xi\in \calD([0,\iy),\R_+):\xi \text{ is nondecreasing}\}.
\]
Endow both spaces with the uniform-on-compacts topology.
Given $\psi^1,\psi^2\in \calP$ and $\zeta,\vr \in \calE$, the {\it dynamics associated
with the initial condition $x$ and the data $\psi,\zeta,\vr$} is given by
\begin{equation}\label{eq022}
\ph(t)= x+yt+ \psi^1(t)-\psi^2(t)+\zeta(t)-\vr(t),\qquad t\ge0.
\end{equation}
We sometimes write the dependence of the dynamics on the data as
$\ph[\psi,\zeta,\vr]$.
Note the analogy between the above equation and equation \eqref{eq011},
and between the condition that $\zeta$ and $\vr$ are nondecreasing and property \eqref{eq012}. The control $\zeta$
stands for the scaled idle time process $\tilde{Z}^n$ and $\vr$ stands for the scaled rejection process $\tilde{R}^n$.
The following condition, analogous to property \eqref{eq013}, will also be in use, namely
\begin{equation}\label{eq023}
\ph(t)\in[0,D],\qquad t\ge 0.
\end{equation}
The game is defined in the sense of Elliott and Kalton \cite{Elli-Kal},
for which we need the notion of strategies.
A measurable mapping $\alpha:\calP^2\to \calE^2$ is called a {\it strategy for the minimizing player}
if it satisfies a causality property. Namely,
for every $\psi=(\psi^1,\psi^2), \tilde \psi=(\tilde\psi^1,\tilde\psi^2) \in \calP^2$
and $t\in[0,\iy)$,
\begin{equation}\label{eq024}
\text{$(\psi^1,\psi^2)(s)=(\tilde\psi^1,\tilde\psi^2)(s)$
for every $s\in[0,t]$\quad implies\quad $\alpha[\psi](s)=\alpha[\tilde\psi](s)$
for every $s\in[0,t]$.}
\end{equation}
Given an initial condition $x$,
a strategy $\al$ is said to be {\it admissible} if, whenever $(\psi^1,\psi^2)\in \calP^2$ and
$(\zeta,\vr)=\al[\psi]$, the corresponding dynamics \eqref{eq022} satisfies the buffer size
constraint \eqref{eq023}. The set of all admissible strategies for the minimizing
player is denoted by $\A$ (or, when the dependence on the initial condition is
important, $\A_x$).
Given $T\in\R_+$, $x\in[0,D]$, $\psi=(\psi^1,\psi^2)\in \calP^2$, and $(\zeta,\vr)\in \calE^2$, we define the cost until time $T$ by
$$
c(x,T,\psi,\zeta,\vr):=\int_0^Th(\ph(t))dt+r\vr(T)-\Ir(T,\psi),
$$
where $\ph$ is the corresponding dynamics.
The value of the game is defined by
$$
V(x):=\inf_{\alpha\in \A_x}\sup_{\psi\in \calP^2, T\in\R_+}c(x,T,\psi, \alpha[\psi]).
$$
We call $\psi$ the \emph{path control} and the $T$ a \emph{time control}, or sometimes
the {\it termination time}. Note that both are controlled by the maximizer.

\subsection{The multi-class model}

In a similar manner to the single-class model, we describe a model with
$I\ge1$ classes of customers, where a buffer of finite length is dedicated
to each class. The various primitives, parameters and processes are denoted with
an additional index $i\in\calI:=\{1,\ldots,I\}$, representing the class. As far as the control
is concerned, an important difference from the single-class model is that
the DM has to determine the allocation of effort to the various classes, on top of
making admission decisions.
Thus, with $\mathbb{S}=\{x=(x_1,\ldots,x_I)\in[0,1]^I:\sum x_i\le1\}$, let
$B^n$ be a process taking values in $\mathbb{S}$, whose $i$-th component represents
the fraction of effort devoted by the server to the class-$i$ customer at the head of the line.
Then the number of service completions of class-$i$ jobs during
the time interval $[0,t]$ is given by $S^n_i(T^n_i(t))$,
where
\begin{equation}\label{eq1007}
T^n_i(t):=\int_0^t B^n_i(s)ds
\end{equation}
is the time devoted to class-$i$ customers by time $t$.
We have the balance equation
\begin{equation}\label{eq1008}
X^n_i(t)=X^n_i(0)+A^n_i(t)-S^n_i(T^n_i(t))-R^n_i(t).
\end{equation}
The process $U^n:=(B^n, R^n)$ is now regarded as a control.
For $i\in\I$ fixed constants $D_i>0$ serve in defining the buffer sizes,
and so the definition of an admissible control now requires that
for every $i\in\I$ one has
\begin{equation}
X^n_i(t)\in[0,b_n\sqrt{n}D_i],\qquad t\geq 0.\notag
\end{equation}
Scaling that is applied in a similar fashion to that in the single-class case
gives rise to the equation
\begin{equation} \label{eq1013}
 \tilde{X}^n_i(t)=\tilde{X}^n_i(0) + y^n_it
 +\tilde{A}^n_i(t)-\tilde{S}^n_i(T^n_i(t))
 + Z^n_i(t)- \tilde{R}^n_i(t),
\end{equation}
where we denote $\rho_i=\la_i/\mu_i$ and
\begin{equation}\notag
Z^n_i(t):=\frac{\mu^n_i}{n}\frac{\sqrt{n}}{b_n}(\rho_i t-T^n_i(t)),
\qquad y^n_i:=\tilde\la^n_i-\rho_i\tilde\mu^n_i.
\end{equation}
The critical load condition now takes the form $\sum_i\rho_i=1$.

An important role is played by the property
\begin{equation}\label{eq1015}
\sum_i\frac{n}{\mu^n_i}Z^n_i \quad \text{starts from zero and is nondecreasing,}
\end{equation}
that holds thanks to the fact that $\sum_iB^n_i\le1$ while $\sum_i\rho_i=1$.
The probabilistic assumptions on the primitives are analogous to those previously
imposed, where in addition, the $2I$ primitive processes are assumed mutually independent.
Finally, the cost is assumed to be of the form
\begin{equation}\notag
J^n(U^n,T):=\frac{1}{b_n^2}
\log\E\Big[\int_0^T e^{b^2_n[\int_0^t \uh\cdot\tilde{X}^n(s))ds+\ur\cdot\tilde{R}^n(t))]}dt\Big],
\end{equation}
where $\uh$ and $\ur$ are given positive vectors and the value is given by
\begin{equation}\notag
V^n(T)=\inf_{U^n} J^n(U^n,T).
\end{equation}

Without going into details, we present the form of the DG that is obtained from
the balance equation \eqref{eq1013}, the condition \eqref{eq1015} and the cost structure;
the derivation is very similar to that presented earlier.
Let $\utheta=(\frac{1}{\mu_1},\ldots,\frac{1}{\mu_I})$ and $\uy=(y_1,\ldots,y_I)$
where $y_i=\tilde\la_i-\rho_i\tilde\mu_i$.
Denote
\begin{equation}\notag
\calE_{\utheta}=\{\uzeta\in \calD([0,T],\R_+^I):\utheta\cdot\uzeta \text{ is nondecreasing}\}.
\end{equation}
Given $\upsi=(\upsi^1,\upsi^2)\in \calP^{2I}$, $\uzeta\in \calE_{\utheta}$, and $\uvr\in\calE^I$, the {\it dynamics, $\uph=(\ph_1,\ldots,\ph_I)$, associated
with the initial state $\ux\in\prod_{i=1}^I[0,D_i]$ and the data $\upsi,\uzeta,$ and $\uvr$} is given by
\begin{equation}\label{eq1028}
\varphi_i(t)= x_i+y_it+ \psi^1_i(t)-\psi^2_i(t)+\zeta_i(t)-\vr_i(t),\qquad t\ge 0,\quad i\in\calI.
\end{equation}
Note the analogy between the above equation and equation \eqref{eq1013},
and between the condition that $\utheta\cdot\uzeta$ is
nondecreasing and property \eqref{eq1015}.
The buffer constraint is now translated into the condition
\begin{equation}\notag
 \varphi_i(t)\in[0,D_i],\qquad t\ge0,\quad i\in\calI.
\end{equation}
A measurable mapping $\ual:\calP^{2I}\to \calE_{\utheta}\times\calE^I$ is called a {\it strategy for the minimizing player}
if it satisfies a causality property. Namely,
for every $\upsi, \tilde \upsi \in \calP^{2I}$
and $t\in[0,T]$,
\begin{equation}\notag
\text{$\upsi(s)=\tilde\upsi(s)$
for all $s\in[0,t]$\quad implies\quad $\ual[\upsi](s)=\ual[\tilde\upsi](s)$
for all $s\in[0,t]$.}
\end{equation}
As before, we denote by $\A_{\ux}$ the set of all admissible strategies given the initial state $\ux$. The cost is given by
\begin{equation}\notag
\uc(\ux,T,\upsi,\uzeta,\uvr):=\int_0^T\uh\cdot \uph(t)dt+r\cdot\uvr(T)-\Jr(T,\upsi),
\end{equation}
where $\uph$ is the corresponding dynamics given in \eqref{eq1028} and $\Jr$ the corresponding rate function (specified in \eqref{100} in the appendix).
The value of the game is defined by
\begin{equation}\notag
V_{MD}(x):=\inf_{\ual\in \A_{\ux}}\sup_{\upsi\in \calP^{2I}, T\in\R_+}c(\ux,T,\upsi, \ual[\upsi]).
\end{equation}
Our result concerning this DG is that {\it it can be reduced to the one introduced earlier}
for the single-class model. In the heavy traffic literature, this type of result is often called
{\it state space collapse}. The precise statement of this reduction, along with
its proof are provided in the appendix. Since this multidimensional version
of the game can be reduced to a one-dimensional one, the rest of this paper
is concerned with the latter.

\section{Solution and analysis of the game}\label{sec_3}
\beginsec

In this section we provide a solution of the game.
After proving some basic properties in Section \ref{sec_3a},
we present, in Section \ref{sec_3b}, an equivalent {\it hitting time} game which is easier to analyze.
In Section \ref{sec_3c} we present the Bellman equation associated with the game.
In Section \ref{sec_3e} we characterize the value function as a viscosity solution of the Bellman equation.
An explicit expression for the minimizer's optimal strategy
appears in Section \ref{sec_3f}.
Finally, a saddle point result appears in Section \ref{sec_3d}.

The minimizer's strategy that is shown to be optimal is of a $\beta$-barrier form.
Informally, this is a strategy that uses the minimal control $(\zeta,\vr)$
so as to keep the dynamics $\ph$ in $[0,\beta]$ at all times. The definition is
based on the {\it Skorohod map on an interval}.
To introduce this map, fix $a<b$. The Skohorod map on the interval $[a,b]$,
denoted by $\Gam_{[a,b]}$, is a map $\calD([0,\iy),\R)\to \calD([0,\iy),\R^3)$.
It is characterized as the solution
map $\om\to(\ph,\eta_1,\eta_2)$ to the so called {\it Skorohod Problem},
namely the problem of finding, for a given $\om$, a triplet $(\ph,\eta_1,\eta_2)$, such that
\begin{equation}\label{44}
\ph=\om+\eta_1-\eta_2,\qquad \ph(t)\in[a,b] \text{ for all } t,
\end{equation}
$\eta_i$ are nonnegative and nondecreasing, $\eta_i(0-)=0$, and
\begin{equation}\label{45}
\int_{[0,\iy)}1_{(a,b]}(\ph)d\eta_1=\int_{[0,\iy)}1_{[a,b)}(\ph)d\eta_2=0.
\end{equation}
By writing $\eta_i(0-)=0$ we adopt the convention that $\eta_i(0)>0$ is regarded a jump
at zero. This convention, in conjunction with $\int_{[0,\iy)}1_{(a,b]}(\ph)d\eta_1=0$
(resp., $\int_{[0,\iy)}1_{[a,b)}(\ph)d\eta_2=0$),
means that if $\om(0)<a$ (resp., $\om(0)>b$) then $\ph(0)=a$ (resp., $b$).
If, however,
$\om(0)\in[a,b]$ then $\ph(0)=\om(0)$, and $\eta_i$ have no jump at zero.

See \cite{KLRS} for existence and uniqueness of solutions,
and continuity and further properties of the map.
In particular, it is well-known that $\Gam_{[a,b]}$ is continuous in the uniform-on-compacts
topology.

\begin{definition}
Fix $(x,\beta)\in[0,D]^2$. The strategy $\al_\beta=(\al_{\beta,1},\al_{\beta,2})\in\A_x$
is called a $\beta$-barrier strategy if for every $\psi\in\calP^2$ one has
$(\ph,\al_{\beta,1},\al_{\beta,2})[\psi]=\Gam_{[0,\beta]}(\psi)$.
\end{definition}

The specific value $\beta_0$ of $\beta$ that gives rise to an optimal strategy
will be determined later, based on the free boundary associated with the Bellman equation.
For the moment we only comment that it is possible to select $\beta=0$, and in fact this choice
turns out to be useful in Section \ref{sec_3a} below.

\subsection{Basic properties}\label{sec_3a}

We begin by providing a necessary and sufficient condition
for the value function to be finite and presenting further properties that are used in the sequel.
\begin{lemma}\label{lem_V_iy}
i. If $-y<r/(4\cc)$ then for every $x\in[0,D]$ one has $V(x)=\iy$.\\
\noi ii. If $-y\ge r/(4\cc)$ then for every $x\in[0,D]$ one has $0\leq V(x)\leq rx$, in particular $V(0)=0$. Moreover, $V$ is nondecreasing. Finally, $V$ is continuous.
\end{lemma}

\noi{\bf Proof of Lemma \ref{lem_V_iy}:}
i. Fix $x\in[0,D]$. Consider a maximizer's control given by $T^\sharp\in\R_+$ and
$\psi^{\sharp,1}(t)=rt/(2\cco)$, $\psi^{\sharp,2}(t)=-rt/(2\cct)$, $t\ge0$. Let $\al\in\A_x$, and denote by $(\zeta^\sharp,\vr^\sharp)=\al[\psi^\sharp]$. The dynamics under $(\psi^{\sharp,1},\psi^{\sharp,2})$ and $(\zeta^\sharp,\vr^\sharp)$ is denoted by $\ph^\sharp$.
We show that the cost is bounded below by a quantity that converges to infinity as $T^\sharp\rightarrow\iy$. Since the maximizer is free to choose $T^\sharp\in\R_+$, the supremum of the cost function over $T$ is infinity. Indeed,
\begin{align}\label{eq026}
c(x,T^\sharp,\psi^\sharp,\al^\sharp[\psi^\sharp])
&=
\int_0^{T^\sharp} h(\ph^\sharp(t))dt+r\vr^\sharp(T)-\Ir(T^\sharp,\psi^\sharp)\\\label{eq027}
&\ge
r\vr^\sharp(T^\sharp)-\Ir(T^\sharp,\psi^\sharp)
=
r\vr^\sharp(T^\sharp)- \frac{1}{4\cc}r^2T^\sharp\\\label{eq028}
&\ge
r(x+\frac{1}{2\cc}rT^\sharp+yT^\sharp+\zeta^\sharp(T)-D )- \frac{1}{4\cc}r^2T^\sharp\\\label{eq029}
&=
r(x-D) + (\frac{1}{4\cc}r^2+yr)T^\sharp+r\zeta(T^\sharp).
\end{align}
Inequality \eqref{eq027} follows since $h$ is nonnegative. Inequality (\ref{eq028}) follows by the dynamics of the game (see \eqref{eq022}) and constraint \eqref{eq023}. Equation \eqref{eq029} is merely a rearrangement of the terms.

\skp

\noi
ii. Fix $x\in[0,D]$. To see why $0\leq V(x)\leq rx$ we show that the maximizer (resp., minimizer) can force a cost of at least $0$ (resp., at most $rx$).
The maximizer may choose $T=0$, in which case for any strategy $\al\in\A_x$, the cost is bounded
below by $0$. This shows $V(x)\ge0$.
The minimizer can choose a $0$-barrier strategy so that the dynamics remain at zero.
That is, the minimizer may choose a strategy $\al_0$
that, given $T\in\R_+$ and $\psi\in\calP^2$ responds with $(\zeta_0,\vr_0)$ where
\begin{align}\notag
\dot\zeta_0(t) &= (\dot\psi^1(t)-\dot\psi^2(t)+y)^-,\qquad 0< t\le T,\\\notag
\dot\vr_0(t) &=  (\dot\psi^1(t)-\dot\psi^2(t)+y)^+,\qquad 0< t\le T,
\end{align}
with $\zeta_0(0)=0$ and $\vr_0(0)=x$.
Fix $T\in\R_+$ and $\psi\in\calP^2$. From the definition of $(\zeta_0,\vr_0)$, and since for every $t>0$ one has $h(\ph(t))=h(0)=0$ it follows that
\begin{align}\label{eq032}
c(x,T,\psi,\al_0[\psi])&=\int_0^Th(\ph(t))dt+r\vr_0(T) -\Ir(T,\psi) \\\notag
&=r\vr_0(0)+ \int_0^T[h(\ph(t))+r(\dot\psi^1(t)-\dot\psi^2(t)+y)^+ \\\notag
&\qquad\qquad\qquad\quad-\cco(\dot\psi^1)^2(t)-\cct(\dot\psi^2)^2(t)]dt\\\notag
&=r\vr_0(0)+\int_0^T[r(\dot\psi^1(t)-\dot\psi^2(t)+y)^+
-\cco(\dot\psi^1)^2(t)-\cct(\dot\psi^2)^2(t)]dt.
\end{align}
One can easily verify that for every $\dot\psi^1(t),\dot\psi^2(t)\in\R$ the expression
given in the integrand is bounded above by $r^2/(4\cc)+ry$, which is nonpositive by assumption.
Therefore the cost given in \eqref{eq032} is bounded above by $r\vr(0)=rx$.
This shows $V(x)\le rx$.

From the inequality $0\leq V(x)\leq rx$ applied to $x=0$ it follows that $V(0)=0$. 

To see why $V$ is nondecreasing fix $0\leq x_1<x_2\leq D$. Given any strategy $\al_2$ associated with the initial state $x_2$, the minimizer may choose for the initial state $x_1$
a strategy $\al_1$ obtained from $\al_2$ by adding to the $\zeta$ component the constant
$x_2-x_1$ (without incurring any additional cost), leaving the $\vr$ component
as under $\al_2$. Therefore, the optimal cost associated with the initial state $x_1$ is less than or equal to that associated with $x_2$. That is, $V(x_1)\leq V(x_2)$.

Since $V$ is nondecreasing its continuity follows by showing that it has no upwards jumps. Assume to the contrary that $V$ has a jump of size $M>0$ at $x\in [0,D]$. In case, $x\in (0,D)$ then for every sufficiently small $\delta>0$ one has
$V(x-\delta)\le V(x+\delta)-M$. W.l.o.g.\ assume that $\delta <\min\{x, M/(4r)\}$.
Therefore, by the definition of $V$ applied to $x-\delta$ it follows that there is a strategy $\al=(\al_1,\al_2)\in\A_{x-\delta}$ such that
\begin{align}\label{eq032aa}
\sup_{\psi\in \calP^2, T\in\R_+}c(x-\delta,T,\psi, \alpha[\psi])<V(x-\delta)+M/2.
\end{align}
Define a new strategy $\hat\al=(\hat\al_1,\hat\al_2)\in\A_{x+\delta}$ by $(\hat\al_1,\hat\al_2)=(\al_1,\al_2-2\delta)$ then clearly,
\begin{align}\label{eq032ab}
\sup_{\psi\in \calP^2, T\in\R_+}c(x+\delta,T,\psi, \hat\alpha[\psi]) =2\delta+ \sup_{\psi\in \calP^2, T\in\R_+}c(x-\delta,T,\psi, \alpha[\psi]).
\end{align}
From \eqref{eq032aa} and \eqref{eq032ab} and the choice of $\delta$ it follows that
$\sup_{\psi\in \calP^2, T\in\R_+}c(x+\delta,T,\psi, \hat\alpha[\psi]) <V(x+\delta)$, which contradicts the definition of $V$. In case $x=0$ or $x=D$ the proof requires only minor modifications and therefore omitted.

\noi\hfill $\Box$

\noi In the rest of the paper,
our standing assumption will be that $-y\ge r/(4\cc)$ holds; thus
by Lemma \ref{lem_V_iy} the value function is finite.

\subsection{The hitting time game}\label{sec_3b}

In the proof of Lemma \ref{lem_V_iy} we saw that if the initial state of the dynamics, $\ph(0)=x$, equals zero then the maximizer prefers to terminate the game immediately. Therefore, it seems plausible that, for a general initial condition, if both players play in an optimal way then as soon as the dynamics hits zero, the game is terminated. Hence the minimizer has no interest in ever activating the $\zeta$ component, since doing so only increases the running cost $h(\ph(t))$ and postpones the termination of the game. In this section we validate this intuition
by formulating an equivalent game that admits the following modifications: first, the $\zeta$ component of the minimizer's control vanishes. Second, the game terminates
when the dynamics hit zero.

Given $\psi\in\calP^2$ and $\vr\in\calE$, the dynamics of the new game are given by
$$
\ph =\ph_{x,\psi,\vr}(t):=x+yt+\psi^1(t)-\psi^2(t)-\vr(t),\qquad t\ge0.
$$
A strategy is a map $\bar\al:\calP^2\rightarrow \calE$ that satisfies the following version of the causality property:
$$
\text{$\psi(s)=\tilde\psi(s)$\ for every $s\in[0,t]$
\ implies\  $\bar\al[\psi](s)=\bar\al[\tilde\psi](s)$
for every $s\in[0,t]$.}
$$
A strategy is said to be admissible if for every $t\ge0$ the dynamics satisfies $\ph(t)\in[0,D]$.
Denote the set of all such strategies by $\bar \A_x$.
Let
\begin{align}\label{eq037}
\tau=\tau(x,\psi,\vr):= \inf\{t\ge0 : \ph_{x,\psi,\vr}(t)=0\},
\end{align}
where throughout we use the convention $\inf\emptyset=\iy$.
For $x\in[0,D]$ and $\bar\al\in\bar\A_x$, define
$$\calQ[x,\bar\al]:=\{\psi\in\calP^2 : \tau(x,\psi,\bar\al[\psi])<\iy\}.$$
This is the set of all the controls for which the dynamics of the game hits zero in a finite time.
Given $x\in[0,D]$, $\bar\al\in\bar\A_x$, and $\psi\in \calQ[x,\bar\al]$, we define the cost by
$$
\bar c(x,\psi,\bar\al[\psi]):=\int_0^\tau h(\ph(t))dt+r\bar\al[\psi](\tau)-\Ir(\tau,\psi),
$$
where $\tau=\tau(x,\psi,\bar\al[\psi])$.
The value of the hitting time game is defined by
$$
U(x)=\inf_{\bar\al\in \bar \A_x}\sup_{ \psi\in \calQ[x,\bar\al]}\bar c(x,\psi, \bar\al[\psi]).
$$
The following proposition states that the value functions of both games are equal to each other.
\begin{proposition}\label{prop_equiv_game}
If $-y\ge r/(4\cc)$ holds then for every $x\in[0,D]$ one has $V(x)=U(x)$.
\end{proposition}
\noi{\bf Proof:}
For every $x\in[0,D]$, every $\psi=(\psi^1,\psi^2)\in\calP^2$, and every $\al=(\al_1,\al_2)\in\A_x$ let
$$\pi=\pi(x,\psi,\al):=\inf\{t\ge0 : x+yt+\psi^1(t)-\psi^2(t) +\al_1[\psi](t)-\al_2[\psi](t)=0\}.$$ 
This is the first time that the dynamics in the original game hit zero. Let $\A_{x,0}$ be the set of all admissible strategies $\al:\calP^2\to\calE^2$ such that $\al^1[\psi]=0$ on $[0,\pi)$ and thereafter $\al$ acts as a $0$-barrier strategy.

We prove the proposition in two steps. First, we show that the infimum in the definition of $V$ can be taken over the set $\A_{x,0}$. That is,
\begin{align}\label{eq041}
V(x)=V_0(x):=\inf_{\al\in\A_{x,0}}\sup_{\psi\in\calP^2, T\in\R_+} c(x,T,\psi,\al[\psi]).
\end{align}
On the second step we show that we can replace the set of controls $\calP^2\times\R_+$ for the maximizer by $\calQ[x,\al_2]$. That is,
\begin{align}\label{eq042}
V_0(x)=V_1(x):=\inf_{\al\in\A_{x,0}}\sup_{\psi\in\calQ[x,\al_2]} \bar c(x,\psi,\al_2[\psi])  .
\end{align}
Once this is achieved, the result follows because it is obvious that $V_1=U$.

\skp\noi \textbf{Step 1:} Fix $x\in[0,D]$.
Since $\A_{x,0}\subseteq\A_x$ it is clear that $V(x)\le V_0(x)$. %

To prove that $V(x)\ge V_0(x)$ we show that for every strategy $\al\in\A_x$ there
exists a strategy $\al^*\in\A_{x,0}$ such that for every $\psi\in\calP^2$ one has
\begin{align}\label{eq043}
\sup_{T\in\R_+} c(x,T,\psi,\al[\psi])\geq \sup_{T\in\R_+} c(x,T,\psi,\al^*[\psi]).
\end{align}
By taking $\sup_{\psi\in\calP^2}$ and then $\inf_{\al\in\A_x}$ it then follows that $ V(x)\ge V_0(x).$
To show \eqref{eq043},
fix a strategy $\al=(\al_1,\al_2)\in\A_x$. Given $\psi\in\calP^2$, let
$$\pi^* = \pi^*(x,\psi,\al):=\inf\{t\ge 0 : x+yt+\psi^1(t)-\psi^2(t)-\al_2[\psi](t)\le 0\}.$$ 
Notice that we omitted $\al_1[\psi]$ from the dynamics. Define $\al^*[\psi](t)=(0,\al_2[\psi](t))$, $0\le t<\pi^*$, and thereafter let $\al^*$ act as a $0$-barrier strategy. One has that $\pi^*(x,\psi,\al)=\tau(x,\psi,\al_2[\psi])$.
By construction, $\al^*_2[\psi](t)\le\al_2[\psi](t)$ for every $t\in[0,\pi^*]$.
Let $\ph$ (resp., $\ph^*$) denote
the dynamics associated with $\psi$ and $\al$ (resp., $\psi$ and $\al^*$). Obviously,
\begin{align}\label{eq044}
\sup_{T\in\R_+}c(x,T,\psi,\al[\psi])&= \sup_{T\in\R_+}\Big\{\int_0^T h(\ph(t))dt +r\al_2[\psi](\pi^*)-\Ir(T,\psi)\Big\}\\\notag
&\geq \sup_{T\in\R_+}\Big\{\int_0^{T\wedge\pi^*} h(\ph(t))dt +r\al_2[\psi](T\wedge\pi^*)-\Ir(T\wedge\pi^*,\psi)\Big\}.
\end{align}
Since over the time interval $[0,\pi^*]$ one has $\al^*_2[\psi]\le \al_2[\psi]$
and $\ph^*\le\ph$, and since $h$ is increasing, it follows that the r.h.s.~of \eqref{eq044}
will not increase if we replace $(\al,\ph)$ by $(\al^*,\ph^*)$.
Therefore, for every $\psi\in\calP^2$ one has $$ \sup_{T\in\R_+} c(x,T,\psi,\al[\psi])\geq \sup_{T\in\R_+} c(x,T\wedge\pi^*,\psi,\al^*[\psi]).$$
Since starting at time
$\pi^*$ the minimizer uses $0$-barrier strategy, it follows by the same arguments as in part (ii) of the proof of Lemma \ref{lem_V_iy} that the maximizer
terminates the game at time $\pi^*$ or earlier. Inequality \eqref{eq043} follows.

\skp\noi\textbf{Step 2:} Recall that $-y\ge r/(4\cc)$ and therefore $V(x)<\iy$. Therefore, in order to show that \eqref{eq042} holds,
it suffices to show that for every $\al\in\A_{x,0}$ for which
$\sup_{\psi\in\calP^2} c(x,T,\psi,\al[\psi])<\iy$, one has
\begin{align}\label{eq047}
\sup_{\psi\in\calP^2,T\in\R_+} c(x,T,\psi,\al[\psi]) = \sup_{\psi\in\calQ[x,\al_2]} \bar c(x,\psi,\al_2[\psi])  .
\end{align}
We do that by showing two sided inequalities. The inequality $\ge$ is obvious since
\begin{align} \label{eq047a}
\sup_{\psi\in\calQ[x,\al_2]} \bar c(x,\psi,\al_2[\psi]) =  \sup_{\psi\in\calQ[x,\al_2]}  c(x,\pi,\psi,\al[\psi])
\end{align}
and
the supremum on the r.h.s.~of \eqref{eq047a} is taken over a smaller set than the supremum on the l.h.s.~of \eqref{eq047}.

To prove the reverse inequality we show that for every $\psi\in\calP^2$ and every $T\in\R_+$ there is $\psi^*\in\calQ[x,\al_2]$ such that
\begin{align}\label{eq048}
c(x,T,\psi,\al[\psi]) \le  \bar c(x,\psi^*,\al_2[\psi^*]).
\end{align}
Let $\pi=\pi(x,\psi,\al)$ be the first time the dynamics under $\psi$ and $\al$ hit zero, possibly infinity. We argue separately for $T$ that is greater than $\pi$ and for $T$ that is smaller than or equal to $\pi$. For $T>\pi$, let $\psi^*=\psi$. By the definition of $\A_{x,0}$ $\pi(x,\psi,\al)=\tau(x,\psi^*,\al_2[\psi^*])$. Denote this time by $\tau$. The cost $c(x,T,\psi,\al[\psi])$ can be expressed as the sum of the cost until time $\tau$ and the cost between times $\tau$ and $T$.
The first cost is $c(x,\tau,\psi,\al[\psi])=\bar c(x,\psi^*,\al_2[\psi^*])$, and the second cost is non positive since the minimizer uses a $0$-barrier strategy after time $\tau$ and in the proof of Lemma \ref{lem_V_iy} we showed that the maximizer optimal cost under this strategy is zero.

For $T\le\tau$ define $\psi^*=\psi(\cdot\w T)$.
The dynamics associated with $\psi^*=(\psi^{*,1},\psi^{*,2})$ and $\al$ until time $\tau^*=\tau(x,\psi^*,\al_2)$ is $\ph^*(t)=x+yt+\psi^{*,1}(t)-\psi^{*,2}(t)-\al_2[\psi](t)$. Since $y<0$ and $\psi^*$ is constant on the time interval $[T,\iy)$ it follows that $\ph^*$ hits zero in a finite time, and therefore $\psi^*\in\calQ[x,\al_2]$.
Thus
\begin{align}\notag
\bar c(x,\psi^*,\al_2[\psi^*]) &= c(x,\tau^*,\psi^*,\al_2[\psi^*])\\\notag
 &= c(x,T,\psi^*,\al_2[\psi^*]) + \int_T^{\tau^*} h(\ph^*(t))dt+r(\al_2[\psi^*](\tau^*)-\al_2[\psi^*](T))\\ \notag
   &\qquad - (\Ir(\tau^*, \psi^*)-\Ir(T, \psi^*))\\\notag
&\le c(x,T,\psi^*,\al_2[\psi^*]) = c(x,T,\psi,\al_2[\psi^*]).
\end{align}
The second equality expresses merely the sum of the cost until time $T$ and the cost between times $T$ and $\tau^*$. The inequality follows since $\al_2$ and $h$ are nondecreasing and since $\dot\psi^1=\dot\psi^2=0$ on the time interval $[T,\tau^*]$. The last equality follows since on the time interval $[0,T]$ we chose $\psi^*$ to agree with $\psi$.

\noi
\hfill $\Box$

\subsection{Derivation of the Bellman equation}\label{sec_3c}
Here we give an intuitive approach towards the solution of the game, present the
corresponding Bellman equation and discuss its solvability.
The derivation of the Bellman equation uses the hitting time game.

We first write the chain rule under the dynamics
\[
\ph(t)=x+yt+\psi^1(t)-\psi^2(t)-\vr(t),\qquad 0\le t\le T,
\]
where $x\in[0,D]$, $T\in\R_+$, $\psi\in\calP^2$, and $\vr\in\calE$. Namely, for any function $f\in\calC^1([0,D],\R)$ one has
\begin{align}
f(\ph(T))=f(x) + \int_0^T f'(\ph(t))(d\psi^1(t)-d\psi^2(t) + ydt-d\vr^c(t)) + \sum_{0\le t\le T}\Delta f(\ph)(t),
\end{align}
where $\vr^c$ is the continuous part of $\vr$. Thus
\begin{align}\label{eq052}
&\int_0^T\left[ - \cco(\dot\psi^1)^2(t)-\cct(\dot\psi^2)^2(t)+h(\ph(t))\right]dt  +r\vr(T)+f(\ph(T)) \\\notag
&=
f(x) +  \int_0^T \left[- \cco(\dot\psi^1)^2(t)-\cct(\dot\psi^2)^2(t)+h(\ph(t)) + f'(\ph(t))(\dot\psi^1(t)-\dot\psi^2(t) + y)\right]dt  \\\notag
&\quad+\int_0^T (r-f'(\ph(t)))d\vr(t) .
\end{align}
Assuming that there exists
$\bar\al\in\bar\A_x$ that is an optimal strategy, by the definition of $U$
one has
\begin{align}\label{eq053}
\sup_{\psi\in\calQ[x,\bar\al]}\left[\int_0^\tau[ - \cco(\dot\psi^1)^2(t)-\cct(\dot\psi^2)^2(t)+h(\ph(t))]dt  +r\bar\al[\psi](\tau)\right] =U(x).
\end{align}
Noting that $U(\ph(\tau))=U(0)=0$, we have by equations \eqref{eq053} and \eqref{eq052} applied to $f=U$ and $T=\tau$,
\begin{align}\label{eq054}
&\sup_{\psi\in\calQ[x,\bar\al]}\Big\{\int_0^T \left[- \cco(\dot\psi^1)^2(t)-\cct(\dot\psi^2)^2(t)+h(\ph(t))
+U'(\ph(t))(\dot\psi^1(t)-\dot\psi^2(t) + y)\right]dt
\\ \notag
&\hspace{11em} +\int_0^T (r-U'(\ph(t)))d\bar\al[\psi](t)
\Big\}=0.
\end{align}
This gives rise to the following equation in the state variable, namely
\begin{equation}\label{eq055}
\max_{p,q\in\R}\min \Big\{ - \cco p^2 -\cct q^2 +h(x) + f'(x)(p-q + y), r-f'(x)  \Big\}=0.
\end{equation}
The maximum above is achieved by $(p,q)=(rf'(x)/(2\cco),rf'(x)/(2\cct))$ and thus we
obtain
\begin{equation}\label{eq056}
\text{[The Bellman equation]}\qquad\max \Big\{ \calL f'(x)  -h(x), \calH f'(x)  \Big\}=0,
\end{equation}
where $\calL$ and $\calH$ denote the operators
\begin{align}\notag
\calL p=- \frac{1}{4\cc}p^2 -yp,
\qquad\calH p=p - r,\qquad p\in\R.
\end{align}

As suggested by the way that \eqref{eq056} is derived from the dynamic description \eqref{eq054},
rejections should take place at times when the state belongs to the part of the space
where $\calH U'=0$. Heuristically, it is plausible that rejections should occur
only when the state exceeds a threshold, and so it is to be expected that first
and second arguments on the l.h.s.\ of equation \eqref{eq056} vanish for $x<\beta$
and $x>\beta$, respectively, where $\beta$ is the free boundary point.
This also suggests that the minimizer plays optimally
by selecting a $\beta$-barrier strategy.

To get information on the free boundary point
of \eqref{eq056}, note that $\calL U'-h=\calH U'=0$ has to hold at this point (at least
in case when the solution is $\calC^1$).
We thus let the free boundary point, $\beta_0$,
be defined as the unique solution $\beta$ of the equation
\begin{equation}\label{eq025}
\frac{r^2}{4\cc}+ry+h(\beta)=0,
\end{equation}
in case it exists.

Note that there exists a unique solution in case $r^2/(4\cc)+ry+h(0)\le 0$ and
$r^2/(4\cc)+ry+h(D)\ge0$, as follows from the assumption that
$h$ is strictly increasing. Otherwise, equation \eqref{eq025} has no solution.
When $r^2/(4\cc)+ry+h(D)< 0$ we define $\Barr=D$.
In case $r^2/(4\cc)+ry+h(0)=r^2/(4\cc)+ry> 0$ Lemma \ref{lem_V_iy}  tells us that the value function is infinite and there is
no need to define $\beta_0$. Thus
\begin{equation}
  \label{02}
  \beta_0=\begin{cases}
    \ds h^{-1}\Big(\frac{-r^2}{4c}-ry\Big), & \ds -h(D)\le \frac{r^2}{4c}+ry\le-h(0),\\ \\
    D, & \ds\frac{r^2}{4c}+ry<-h(D).
  \end{cases}
\end{equation}

As for a solution to the Bellman equation, we let $g:[0,D]\to\R$ be defined by
\[
g(x)=\left\{\begin{array}{ll}
               \displaystyle\int_0^x 2c\Big(-y-\sqrt{y^2-\frac{h(u)}{\cc}}\Big)du,
               &\ 0\leq x\leq \Barr,
\\ \\
                g(\Barr)+r(x-\Barr), & \  \Barr<x\leq D,
              \end{array}
\right.
\]
where in the case $\beta_0=D$, $g$ is defined on $[0,D]$ through the first line alone.
Then the expression under the root is always positive. To see this notice that from the discussion above and the definition of $\Barr$ it follows that in case $-y\ge r/(4\cc)$ one has $r^2/(4\cc)+ry+h(\Barr)\le 0$. By simple algebraic manipulation it follows that $$(r/(2\sqrt{\cc})+y\sqrt{\cc})^2-y^2\cc+h(\Barr) =r^2/(4\cc)+ry+h(\Barr) \le 0,$$ which yields $y^2-h(\Barr)/\cc\ge 0$. By the monotonicity of $h$ it follows that for every $u\in[0,\Barr]$ one has $y^2-h(u)/\cc\ge 0$.

Moreover, $g$ satisfies the boundary condition $g(0)=0$ and it is easy to check that
it is continuously differentiable on $(0,\beta_0)$ and on $(\beta_0,D)$ (and continuous
on all of $[0,D]$) and satisfies
on these two intervals the two parts of the Bellman equation, namely
$\calL g'-h=0$ and $\calH g'=0$, respectively. Thus if $\beta_0=D$ then $g$ solves
the Bellman equation in a classical sense. However, if $\beta_0<D$ it may occur that
$g$ is not differentiable at $\beta_0$. Specifically, when $\beta_0<D$,
$g$ is not differentiable at $\beta_0$ if and only if
\begin{equation}\label{01}
\frac{r}{4\cc} \le -y < \Big(\frac{r}{4\cc}+\frac{h(D)}{r}\Big)\wedge \frac{r}{2\cc}
\end{equation}
Moreover, under \eqref{01}, the left- and right-derivative at $\beta_0$ satisfy $g'_L(\beta_0)
<g'_R(\beta_0)$.
Although $g$ is not always a classical solution to the Bellman equation,
we will prove in the next section that it is a viscosity solution (suitably defined)
to that equation.

\subsection{Viscosity solutions}\label{sec_3e}

In this section we consider the Bellman equation \eqref{eq056} in viscosity sense and characterize
the value function $U$ of the hitting time game as the unique solution of that equation,
with suitable boundary conditions.
The equality $U=g$ is finally obtained when we prove, in the next section, that $g$
also solves the Bellman equation.

For a set $B\subset\R$ denote by $\oo{B}$ its closure.
\begin{definition}\label{def_viscosity}[Viscosity solution]\\
i. For $B=(0,D)$ or $B=(0,D]$, a continuous function $f:[0,D]\rightarrow\R$
is said to be a viscosity supersolution (resp., subsolution)
of \eqref{eq056} on $B$ if for every $x\in B$ and every $\phi\in\calC^1([0,D],\R)$ for which
$f-\phi$ has a global minimum (maximum) on $B$ at $x$ one has
\[
\max \Big\{ \calL \phi'(x)  -h(x), \calH \phi'(x)  \Big\}\ge0\qquad(\le0).
\]
\\
ii. A function $f:[0,D]\to\R$ is said to be a viscosity solution of \eqref{eq056} on
$[0,D]$ if it is a viscosity subsolution of \eqref{eq056} on $(0,D)$, a viscosity
supersolution of \eqref{eq056} on $(0,D]$, and $f(0)=0$.
\end{definition}
The definition accounts for a Dirichlet boundary condition at $0$ as well as a state constraint
boundary condition at $D$, consisting of the requirement that $f$ is a supersolution up to the boundary
point $D$ (see Section IV.5 of \cite{bardi}).

\begin{theorem}\label{thm_viscosity_V}
The value function $U$ is the unique viscosity solution of \eqref{eq056} on $[0,D]$.
\end{theorem}

The proof is performed in two steps: we first show that $U$ satisfies the Bellman equation
and then prove uniqueness.

\noi{\bf Proof that $U$ is a viscosity solution:}
By Lemma \ref{lem_V_iy} and the equality $U=V$ established in Proposition \ref{prop_equiv_game},
$U$ is continuous and null at zero. It remains to prove
that $U$ is a subsolution on $(0,D)$ and a supersolution on $(0,D]$.

\skp\noi{\bf Subsolution on $(0,D)$:}
Fix $x\in(0,D)$ and let $\phi\in\calC^1([0,D],\R)$ be such that $U-\phi$ has a global maximum at $x$. We can assume w.l.o.g.~that $U(x)=\phi(x)$. Thus $U\leq\phi$ on $[0,D]$. We need to show that
\begin{align}\label{eq099}
 \calL \phi'(x)  -h (x)\le0
\end{align}
and
\begin{align}\label{eq100}
 \calH \phi'(x)\le0.
\end{align}
Consider $\delta>0$ such that $x-\delta \in[0,D]$. The proof of \eqref{eq100} is based on the inequality $U(x)-U(x-\delta)\leq r\delta$, that is an extension of the two-sided inequality
stated in Lemma \ref{lem_V_iy}.
Let us show this inequality.
Given a strategy $\bar\al\in\bar \A_{x-\del}$ let $\tilde\al\in\bar \A_x$ be defined
as $\tilde\al[\psi]=\del+\bar\al[\psi]$.
Then the responses of the two strategies to $\psi$
are related to each other by $\tilde\ph=\bar\ph$ and $\tilde\vr=\bar\vr+\del$.
As a result,
\[
\bar c(x,\psi,\tilde\al[\psi])=\bar c(x-\del,\psi,\bar\al[\psi])+r\del.
\]
Hence $U(x)\le\sup_\psi\bar c(x-\del,\psi,\bar\al[\psi])+r\del$ and since
$\bar\al$ is arbitrary, we have $U(x)\le U(x-\del)+r\del$.
Therefore $\phi(x)- \phi(x-\delta)\le U(x)- U(x-\delta)\le r\delta$, establishing
\eqref{eq100}.

To prove \eqref{eq099} we let $\delta>0$ be such that $0< x-\delta<x+\delta<D$, define
\begin{align}\label{eq102}
\taud:=\delta\wedge\inf\{t\ge0 : |\ph(t)-x|\ge\delta\}
\end{align}
and use the dynamic programming principle (DPP), that states
\[
U(x)=\inf_{\bar\al\in\bar \A_x}
\sup_{\psi\in\calQ[x,\bar\al^\delta]}
\left[\int_0^{\taud}(-\cco(\dot\psi^1)^2(t)-\cct(\dot\psi^2)^2(t)+h(\ph(t)))dt +r\bar\al[\psi](\taud)
+ U(\ph(\taud))\right].
\]
Consider a strategy $\bar\al$ with the property that for every
$\psi\in\calQ[x,\bar\al^\delta]$ and every $t\in[0,\taud]$ one has
$\bar\al[\psi](t)=0$. By the chain rule,
\begin{align}\notag
U(x)=\phi(x)= - \int_0^{\taud} \phi'(\ph(t))(\dot\psi^1(t)-\dot\psi^2(t)+y)dt + \phi(\ph(\taud)).
\end{align}
Using the DPP, the inequality $U\le\phi$, and the above display, we obtain
\begin{align}\notag
U(x)\leq& \sup_{\psi\in\calQ[x,\bar\al^\delta]}\left[\int_0^{\taud}(-\cco(\dot\psi^1)^2(t)-\cct(\dot\psi^2)^2(t)+h(\ph(t)))dt + U(\ph(\taud))\right]\\\notag
\leq&\sup_{\psi\in\calQ[x,\bar\al^\delta]}\left[\int_0^{\taud}(-\cco(\dot\psi^1)^2(t)-\cct(\dot\psi^2)^2(t)+h(\ph(t)))dt + \phi(\ph(\taud))\right]\\\notag
=& \sup_{\psi\in\calQ[x,\bar\al^\delta]}\Big[\phi(x) + \int_0^{\taud}(-\cco(\dot\psi^1)^2(t)-\cct(\dot\psi^2)^2(t)+h(\ph(t)) \\\notag
&\quad\qquad\qquad\qquad\qquad\quad + \phi'(\ph(t))(\dot\psi^1(t)-\dot\psi^2(t)+y))dt\Big].
\end{align}
Recalling that $U(x)=\phi(x)$, we have
$$\sup_{\psi\in\calQ[x,\bar\al^\delta]} \int_0^{\taud}(-\cco(\dot\psi^1)^2(t)-\cct(\dot\psi^2)^2(t)+h(\ph(t))) + \phi'(\ph(t))(\dot\psi^1(t)-\dot\psi^2(t)+y))dt
\geq 0$$ 
for every $\delta>0$.
Hence
\begin{equation}\label{eq106}
\sup_{\dot\psi^1(0),\dot\psi^2(0)\in\R} \left[-\cco(\dot\psi^1)^2(0)-\cct(\dot\psi^2)^2(0)
+h(x) + \phi'(x)(\dot\psi^1(0)-\dot\psi^2(0)+y)\right]\geq 0 .
\end{equation}
The maximum in \eqref{eq106} is achieved by $(\dot\psi^1(0),\dot\psi^2(0))^*= (\phi'(x)/(2\cco),\phi'(x)/(2\cco))$ and equals
$-\calL\phi'(x) + h(x) $, and \eqref{eq099} follows.

\skp\noi{\bf Supersolution on $(0,D]$:}
Fix $x\in(0,D]$ and let $\phi\in\calC^1([0,D],\R)$ be such that $U-\phi$ has a global minimum at $x$. Assume w.l.o.g.~that $U(x)=\phi(x)$. Thus $U\geq\phi$ on $[0,D]$. We need to show that either
\begin{align}\label{eq107}
 \calL \phi'(x)  -h(x) \ge0
\end{align}
or
\begin{align}\label{eq108a}
 \calH \phi'(x)\ge0.
\end{align}
Arguing by contradiction, assume that neither of the above inequalities hold. Then
one can find $\eps>0$ such that
\begin{align}\label{eq109}
\calH \phi'\le -2\eps \quad\text{and}\quad \calL \phi'  -h \le -2\eps
\end{align}
hold on $\oo B_\del(x):=\{\xi\in [0,D] : |x-\xi|\leq\delta\}$.
Denote $\psi_x(t)=(\phi'(x)t/(2\cco),\phi'(x)t/(2\cco))$, $t\ge 0$.
Recall the definition of $\tau_\del$ from \eqref{eq102}. Let $\bar \A_x(\del)$ denote
the collection of strategies $\bar\al\in \bar \A_x$
that have the property $\ph(\tau_\del)\in\oo B_\del(x)$ for all $\psi$, where
$\ph$ denotes the dynamics under $\psi$ and $\bar\al$.
Pick an arbitrary $\bar\al_x\in\bar\A_x(\del)$
and denote $\vr_x= \bar\al_x[\psi_x]$. The associated dynamics are given by
$\ph_x=x+y\cdot+\psi^1_x-\psi^1_x-\vr_x$. Let $\taud$ be as in \eqref{eq102} with $\ph=\ph_x$.
An application of the chain rule gives
\begin{align}\label{eq110}
&\phi(\ph_x(\taud))-\phi(x)\\\notag
&\quad=-\int_0^{\taud}(-\cco(\dot\psi^1_x)^2(t)-\cct(\dot\psi^2_x)^2(t)+h(\ph_x(t))dt
-r\vr_x(\taud)\\\notag%
&\quad\quad+\int_0^{\taud}(-\cco(\dot\psi^1_x)^2(t)-\cct(\dot\psi^2_x)^2(t)+h(\ph_x(t)) + \phi'(\ph_x(t))(\dot\psi^1_x(t)-\dot\psi^1_x(t)+y))dt\\\notag
&\quad\quad+\int_{[0,\taud]}(r-\phi'(\ph_x(t)))d\vr_x(t).
\end{align}
We will bound from below the last three terms on the r.h.s.~of \eqref{eq110}.
By \eqref{eq109}, the definition of $\psi_x$, and the continuity of the functions $\phi'$ and $h$ it follows that there is sufficiently small $\delta>0$ such that for every $t\in[0,\taud]$ one has
\begin{align}\notag
&-\cco(\dot\psi^1_x)^2(t)-\cct(\dot\psi^2_x)^2(t)+h(\ph_x(t)) + \phi'(\ph_x(t))(\dot\psi^1_x(t)-\dot\psi^1_x(t)+y)\\\notag
&\quad\geq
-\cco(\dot\psi^1_x)^2(0)-\cct(\dot\psi^2_x)^2(0)+h(\ph_x(0)) + \phi'(\ph_x(0))(\dot\psi^1_x(0)-\dot\psi^1_x(0)+y)-\eps\\\notag
&\quad=-\calL\phi'(\ph_x(0))+h(\ph_x(0))-\eps=-\calL\phi'(x)+h(x)-\eps\ge\eps,
\end{align}
where the first equality follows by the definition of $\psi_x$ and by \eqref{eq055}.
Therefore the second term in \eqref{eq110} is bounded below by $\eps\tau_\del$.
By using similar arguments it follows from \eqref{eq109} that for sufficiently small $\delta$
\[
\int_{[0,\taud]}(r-\phi'(\ph_x(t)))d\vr_x(t)\geq \eps\vr_x(\taud).
\]
Using the last two inequalities in \eqref{eq110}, we obtain
\begin{align}\label{eq114}\notag
&\phi(\ph_x(\taud))-\phi(x)\\
&\geq -\int_0^{\taud}(-\cco(\dot\psi^1_x)^2(t)-\cct(\dot\psi^2_x)^2(t)+h(\ph_x(t))dt
-r\vr_x(\tau_\del)
+\eps(\taud+\vr_x(\taud)).
\end{align}
We now show that the sum $\taud+\vr_x(\taud)$ is bounded below by $C(x)\delta$, where $C(x)$ is a positive constant that depends solely on $x$. To this end, assume that  $\taud<\delta$ (otherwise one may simply take $C(x)=1$). Therefore, by the definition of $\bar \A_x(\del)$ $|\ph_x(\taud)-x|=\delta$. Notice that if $\vr_x(\taud)\le\delta/2$ then in case $y+\dot\psi^1_x(t)-\dot\psi^2_x(t)=y+\phi'(x)/(2c)<0$ (resp., $>0$), 
the function $yt +\psi^1_x(t)-\psi^2_x(t)$ must pass $-\del/2$ (resp., $3\del/2$) units over the time interval $[0,\taud]$. Therefore $|yt+\frac{1}{2\cc}\phi'(x)t|$ must pass at least $\del/2$ units. Hence $\taud\geq \{2|\frac{1}{2\cc}\phi'(x)+y|\}^{-1}\delta$. Overall,\footnote{Note that the choice of $C(x)$ is independent of the strategy $\bar\al_x$.}
\begin{align}\notag
\taud+\vr_x(\taud)\geq
\min\left\{\dfrac{\delta}{2},\dfrac{1}{2|\frac{1}{2\cc}\phi'(x)+y|}\delta\right\}
=:C(x)\delta,
\end{align}
and \begin{align}\label{eq116}
\phi(\ph_x(\taud))-\phi(x)\ge
&-\int_0^{\taud}(-\cco(\dot\psi^1_x)^2(t)-\cct(\dot\psi^2_x)^2(t)+h(\ph_x(t))dt
-r\vr_x(\tau_\del)
+\eps C(x)\delta.
\end{align}
By the dynamic programming principle, using the fact $\taud\leq \tau$, it follows that
\begin{align}\notag
U(x)=&\inf_{\bar\al\in\bar\A_x(\del)}\sup_{\psi\in\calQ[x,\bar\al]}\Big[
\int_0^{\taud}[-\cco(\dot\psi^1)^2(t)-\cct(\dot\psi^2)^2(t)+h(\ph(t))]dt
+r\bar\al[\psi](\taud) + U(\ph(\taud))\Big]\\\notag%
\geq&\inf_{\bar\al\in\bar\A_x}\big[
\int_0^{\taud}[-\cco(\dot\psi^1_x)^2(t)-\cct(\dot\psi^2_x)^2(t)+h(\ph_x(t))]dt+r\vr_x(\taud) + U(\ph_x(\taud))\big].
\end{align}
Therefore, there exists a strategy $\al^*\in\bar \A_x(\del)$
such that, denoting $\vr_x^*=\al^*[\psi_x]$,
\begin{align}\label{eq118}
U(x)\geq&\int_0^{\taud}[-\cco(\dot\psi^1_x)^2(t)-\cct(\dot\psi^2_x)^2(t)+h(\ph^*_x(t))]dt
+r\vr^*_x(\taud) + U(\ph^*_x(\taud)) - \frac{1}{2}\eps C(x)\delta,
\end{align}
where $\ph_x^*=x+t\cdot+\psi^1_x-\psi^1_x-\vr^*_x$.
Recalling that $\bar\al_x\in\bar \A_x(\del)$ is arbitrary, letting it be $\al^*$, and
substituting
$\ph_x=\ph^*_x$ in \eqref{eq116}, it follows from inequalities \eqref{eq116}
and \eqref{eq118} that
\[
\phi(\ph^*_x(\taud))-U(\ph^*_x(\taud))\ge\frac{1}{2}\eps C(x)\delta,
\]
which contradicts the fact that $U\ge\phi$.

\skp\noi{\bf Proof of uniqueness of viscosity solutions:}
We prove uniqueness by the comparison principle.
Assume that $v$ is a supersolution on $(0,D]$ and $u$ a subsolution on $(0,D)$.
We show that $u\leq v$.

Arguing by contradiction, let us assume to the contrary, that $u\not\leq v$.
Then there exist
$x\in[0,D]$ and $\delta\in(0,1)$ such that $u(x)-(1+\delta)v(x)>0$.
Fix such $\del$. Fix also
$z\in\arg\max_{[0,D]}\left[u-(1+\delta)v\right]$.
Since $u(0)=v(0)=0$ it follows that $z>0$.

For $\eps>0$ define the following test function
\[
\xi^\eps(x,y)=u(x)-(1+\delta)v(y)- \left(\frac{x-y}{\eps} - \left(-1-z\right)\right)^2 -  \left(y-z\right)^2\quad x,y\in[0,D].
\]
Note that $z + \eps\left(-1-z\right)$ is in $[0,D]$ for every small enough $\eps$.
Also,
\begin{equation}\label{eq122}
\xi^\eps\left(z + \eps\left(-1-z\right),z\right)
\ge
u(z)-(1+\delta)v(z)-mod(C\eps),
\end{equation}
where here and in what follows
$mod$ is the modulus of continuity of $u$ and $v$, and $C$ is the constant $1+D$.
We also have
\begin{equation}\label{eq122b}
\xi^\eps\left(x,y\right)
\le
u(x)-(1+\delta)v(x)+2mod(|x-y|) - \left(\frac{x-y}{\eps} - \left(-1-z\right)\right)^2 -
\left(y-z\right)^2.
\end{equation}
For each $\eps>0$ pick $(\xe,\ye)\in\arg\max_{[0,D]^2}\left[\xi^\eps\left(x,y\right)\right]$. From equations \eqref{eq122} and \eqref{eq122b}, substituting $(x,y)=(\xe,\ye)$
and using $u(z)-(1+\del)v(z)\ge u(\xe)-(1+\del)v(\xe)$, it follows that
\begin{equation}\label{eq124}
\left(\frac{\xe-\ye}{\eps} - \left(-1-z\right)\right)^2 +  \left(\ye-z\right)^2
\le mod(C\eps)+2mod(|\xe-\ye|).
\end{equation}
The above yields two properties that are being used in the sequel.\\
(a) $\lim_{\eps\rightarrow\iy}\xe=\lim_{\eps\rightarrow\iy}\ye=z$.\\
(b) For every $\eps$ small enough one has $\xe<\ye$.

\skp\noi Property (a) follows since the r.h.s.~of equation \eqref{eq124} is bounded above and therefore there is $m>0$ such that $|\xe-\ye|<m\eps$.
This also yields that for sufficiently small $\eps$ the r.h.s.~of \eqref{eq124} is bounded above by $z^2$.
Therefore, $|\xe-\ye- \eps\left(-1-z\right)|\le\eps z$. By these two inequalities it follows that
\begin{align}\notag
\xe
&\le
\ye+\eps\left(-1-z\right)+\eps z
\le
\ye -\eps<\ye,
\end{align}
which proves property (b).

\noi Now consider the maps
\begin{align}\label{eq125}
\phiu(x)&:=(1+\delta)v(\ye)
+ \left(\frac{x-\ye}{\eps}-\left(-1-z\right)\right)^2
+  \left(\ye-z\right)^2,\\\label{eq126}
\phiv(y)&:=\frac{1}{1+\delta}
\left[u(\xe)
- \left(\frac{\xe-y}{\eps}-\left(-1-z\right)\right)^2
-  \left(y-z\right)^2\right].
\end{align}
Then $u-\phiu$ has a maximum at $\xe$ and $v-\phiv$ has a minimum at $\ye$.
Note that by property (a) and the fact $z>0$, $\min(\xe,\ye)>0$ provided $\eps$ is small,
and that by property (b), $\xe<D$. As a result, the definition of
a viscosity subsolution on $(0,D)$ and supersolution on $(0,D]$ yields
\begin{align}\label{eq127}
&\max \Big\{ \calL \phiu'(\xe)  -h(\xe), \calH \phiu'(\xe)  \Big\}\le0\\\label{eq128}
&\max \Big\{ \calL \phiv'(\ye)  -h(\ye), \calH \phiv'(\ye)  \Big\}\ge0.
\end{align}
Now,
\begin{align}\label{eq129}
\phiu'(\xe) &= \frac{2}{\eps}\left(\frac{\xe-\ye}{\eps}-\left(-1-z\right)\right)\\\notag
\phiv'(\ye) &= \frac{1}{1+\delta}\left[\frac{2}{\eps}\left(\frac{\xe-\ye}{\eps}-\left(-1-z\right)\right) - 2\left(\ye-z\right)\right]\\\label{eq130}
&=\frac{1}{1+\delta}\phiu'(\xe) - \frac{2}{1+\delta}\cdot(\ye-z).
\end{align}
By \eqref{eq127} it follows that $\calH\phiu'(\xe)\leq 0$, that is, $\phiu'(\xe)\leq r$. Together with property (a) above and equation \eqref{eq130} it follows that $\calH \phiv'(\ye)<0$. Therefore, \eqref{eq128} yields that

\begin{align}\label{eq129a}
\calL \phiv'(\ye)  -h(\ye)\ge 0
\end{align}
We now claim that the sequence $\{\phiv'(\ye)\}_\eps$ is positive and bounded away from zero. Denote by $q=\lim\phiv'(\ye)=\phiv'(z)$, where property (a) is used. Using again property (a) and \eqref{eq129a} it follows that
$\calL (q)  -h(z)\ge0$. Since $z>0$, we get $h(z)>0$,
and therefore $\calL (q)>0$, which in turn yields $q>0$.
By considering a subsequent we may assume w.l.o.g.~that $\{\phiv'(\ye)\}$ is bounded away from zero.
Hence, for sufficiently small $\eps$ one has
\begin{align}\notag
\Big| \frac{2}{1+\delta}\cdot(\ye-z)\Big| \le \frac{\delta}{2(1+\delta)}\phiv'(\ye),
\end{align}
thanks to the fact that the l.h.s.~converges to zero as $\eps\to 0$ by property (a).
Together with \eqref{eq130} we get that for sufficiently small $\eps$
\begin{align}\label{eq130c}
\phiv'(\ye) \le \frac{1}{1+\delta}\phiu'(\xe) + \frac{\delta}{2(1+\delta)}\phiv'(\ye).
\end{align}
Together with the inequality $\phiu'(\xe)\leq r$ we get that
\begin{equation}\label{eq131}
(1+\dfrac{\delta}{2})\phiv'(\ye)\leq\phiu'(\xe)\leq r.
\end{equation}
If we show that for sufficiently small parameters $\eps$ and $\delta$ one has
\begin{equation}\label{eq132}
\calL (r)  -h(\xe)>0
\end{equation}
and
\begin{equation}\label{eq133}
\calL ((1+\dfrac{\delta}{2})\phiv'(\ye))  -h(\xe)>0
\end{equation}
then by the concavity of $p\mapsto\calL (p)$ and \eqref{eq131}--\eqref{eq133} it follows that $\calL (\phiu'(\xe))  -h(\xe)>0$ in contradiction to \eqref{eq127}.
We now prove that inequalities \eqref{eq132} and \eqref{eq133} hold.

\noi {\bf Proof that inequality \eqref{eq132} holds:}
By property (b) and by the monotonicity of $h$ it follows that $\calL(r)-h(\xe)> \calL(r)-h(\ye)$. Therefore, it is sufficient to show that $\calL(r)-h(\ye)\ge 0$. To this end, assume to the contrary that $\calL(r)-h(\ye)< 0$. Recall that from \eqref{eq131} one has $\phiv'(\ye)<r$. Fix $p_\eps\in(\phiv'(\ye),r)$ sufficiently close to $r$ such that
$\calL(p_\eps)-h(\ye)< 0$. Then
\begin{equation}\label{30}
\max \Big\{ \calL (p_\eps)  -h(y_\eps), \calH (p_\eps)  \Big\}< 0
\end{equation}
Denote $b=\min_{\zeta\in[y_\eps,D]}(v(\zeta)-p_\eps\zeta)$ and consider
the function $\phi(\zeta)=p_\eps\zeta+b$, $\zeta\in[0,D]$.
Pick $z_\eps\in[y_\eps,D]$ for which
$v(z_\eps)=\phi(z_\eps)$. By the definition of $b$ one has $v\ge\phi$ on $[y_\eps,D]$.
Moreover, since $p_\eps>\phi'_v(y_\eps)$ it follows that $v\ge\phi$ in fact holds on
a larger interval, namely on $(y_\eps-c,D]$ for some $c>0$.
Thus we may consider $\phi$ as a test function in the definition
of $v$ as a supersolution on $(0,D]$, by which
\[
\max \Big\{ \calL (p_\eps)  -h(z_\eps), \calH (p_\eps)  \Big\}\ge 0.
\]
This, the fact that $z_\eps\ge y_\eps$, and the monotonicity of $h$
give a contradiction to \eqref{30}.

\noi {\bf Proof that inequality \eqref{eq133} holds:}
From property (a), together with $z>0$ we may assume w.l.o.g.~that for every $\eps$ one has $\ye>0$. Also recall that we assumed that $\phiv'(\ye)>0$. Similarly to \eqref{eq131}, one may take sufficiently small $\eps$ such that
\begin{align}\label{eq134}
&0<\phiv'(\ye)<(1+\frac{\delta}{2})\phiv'(\ye)<(1+\frac{2\delta}{3})\phiv'(\ye)<r
\end{align}
and consequently
\begin{align}\label{eq137}
&0<q<(1+\frac{\delta}{2})q<(1+\frac{2\delta}{3})q\le r,
\end{align}
where the fact that $q>0$ has been proved earlier.
By taking $\eps\rightarrow 0$ in \eqref{eq132} and \eqref{eq129a} it follows form the limits $\lim\xe=\lim\ye=z$ (property (a)) and $\lim\phiv'(\ye)=q$ that
\begin{equation}\label{eq136}
 \calL (q)  -h(z)\ge 0\quad\text{and}\quad\calL (r)  -h(z)\ge 0.
\end{equation}
From \eqref{eq136} and \eqref{eq137} and the strict concavity of $p\mapsto\calL(p)$ it follows that
$
\calL ((1+\dfrac{\delta}{2})q)  -h(z)>0.
$
By the limits $\lim\xe=z$ and $\lim\phiv'(\ye)=q$ it follows that for sufficiently small $\eps$ inequality \eqref{eq133} holds.
\hfill $\Box$

\subsubsection*{Explicit viscosity solution of the Bellman equation}

We now prove that the function $g$ is the value function. We do that by showing that $g$ is a viscosity solution of the Bellman equation, which admits a unique solution.
\begin{theorem}\label{thm_viscosity_g}
The function $g$ is a viscosity solution of \eqref{eq056} on $[0,D]$.
\end{theorem}
Note that this result combined with Theorem \ref{thm_viscosity_V} shows $U=g$.

\skp

\noi
{\bf Proof:}
One verifies directly that the function $g$ satisfies, in the classical sense,
$\calL g'-h=0$ and $\calH g'\le0$
on $(0,\beta_0)$, and (in case $\beta_0<D$) also $\calL g'-h\le 0$ and $\calH g'=0$ on
$(\beta_0,D)$. Also, if $r/(2\cc) \le -y < r/(4\cc)+h(D)/r$ then as follows from Section \ref{sec_3c} $\Barr<D$ and the function $g$ is differentiable at $\Barr$ and one can easily verify that $\calL g'(\Barr)-h(\Barr)= \calH g'(\Barr)=0$. Moreover, $g(0)=0$. Hence it remains to verify the viscosity solution
definition at the points $\beta_0$ and $D$ only.

\skp\noi{\bf Subsolution for $\Barr$:}
The case $r/(2\cc) \le -y < r/(4\cc)+h(D)/r$ was discussed before. We analyze the two remaining cases $ -y \ge r/(4\cc)+h(D)/r$ and $r/(4\cc) \le -y< (r/(4\cc)+h(D)/r)\wedge r/(2\cc)$.
If $ -y \ge r/(4\cc)+h(D)/r$ then as we saw in Section \ref{sec_3c}, $\Barr=D$, and the subsolution property need not be verified since by definition the viscosity subsolution property should be verified only on the interval $(0,D)$.
If $r/(4\cc) \le -y< (r/(4\cc)+h(D)/r)\wedge r/(2\cc)$ then as we saw in Section \ref{sec_3c}, $g_L'(\Barr)<g_R'(\Barr)$. Therefore, there is no function $\phi\in\calC^1([0,D],\R)$ such that $g-\phi$ has a global maximum at $\Barr$, and the claim holds vacuously.

\skp\noi{\bf Supersolution for $\Barr$ and $D$:}
Fix $x\in\{\Barr,D\}$. Let $\phi\in\calC^1([0,D],\R)$ be such that $g-\phi$ has a global minimum at $x$. We can assume w.l.o.g.~that $g(x)=\phi(x)$. Thus, $g\geq\phi$ on $[0,D]$. We need to show that
\begin{align}\label{eq140}
 \calL \phi'(x)  -h(x) \ge0
\end{align}
or
\begin{align}\label{eq141}
 \calH \phi'(x)\ge0.
\end{align}
As before, we analyze the two cases $ -y \ge r/(4\cc)+h(D)/r$ and $r/(4\cc) \le -y< (r/(4\cc)+h(D)/r)\wedge r/(2\cc)$.
If $r/(4\cc) \le -y< (r/(4\cc)+h(D)/r)\wedge r/(2\cc)$ then as was shown in Section \ref{sec_3c} $\Barr<D$, $g_L'(\Barr)<g_R'(\Barr)$, and $g'(D)=r$.
For $x=D$ \eqref{eq141} holds since $\phi'(D)\ge g'(D)$.
For $x=\Barr$. Since
$
g_L'(\Barr)\leq\phi'(\Barr)\leq g_R'(\Barr)
$
it follows that $2\cc(-y-\sqrt{y^2-h(\Barr)/\cc})\leq\phi'(\Barr)\leq r$.
One can verify that
\begin{equation}\label{eq142}
\calL (2\cc(-y-\sqrt{y^2-h(\Barr)/\cc}))  -h(\Barr)  \ge0\quad\text{and}\quad\calL r  -h(\Barr) \ge 0.
\end{equation}
From the above and by the concavity of $\calL$ it follows that \eqref{eq140} holds.

If $ -y \ge r/(4\cc)+h(D)/r$ then as was shown in Section \ref{sec_3c} $\Barr=D$.
If $\phi'(\Barr)\ge r$ then \eqref{eq141} holds. If $\phi'(\Barr)\le r$ then
$2\cc(-y-\sqrt{y^2-h(\Barr)/\cc})=g_L'(\Barr)\leq\phi'(\Barr)\leq r$ and from \eqref{eq142} as before it follows that \eqref{eq140} holds.

\hfill $\Box$

\subsection{Optimal strategy}\label{sec_3f}

The following theorem shows that the $\Barr$-barrier strategy is optimal.
We defined the barrier strategies only for the original game. The definition for the hitting time game is similar.
\begin{definition}
Fix $(x,\beta)\in[0,D]^2$. The strategy $\bar\al_\beta\in\bar\A_x$ is called a $\beta$-barrier strategy if for every $\psi\in\calQ[x,\bar\al_\gamma]$ one has
$(\ph,0,\bar\al_\beta)(\psi)=\Gam_{[0,\gamma]}(\psi)$ on the time interval $[0,\tau]$.
\end{definition}

To state the result let us add a piece of notation.
Denote by $\bar J(x,\bar\al)$ the optimal cost achieved by the maximizer in the hitting time game when the initial state is $x$ and the minimizer uses the strategy $\bar\al\in\bar\A_x$. That is,
\begin{equation}\label{eq068}
\bar J(x,\bar\al):=\sup_{\psi\in\calQ[x,\bar\al]}\bar c(x, \psi,\bar\al[\psi]), \qquad x\in[0,D],\quad \bar\al\in\A_x.
\end{equation}

\begin{theorem}\label{thm_g>J}
Let $\bar\al_\Barr$ be the $\Barr$-barrier strategy. For every $x\in[0,D]$ one has $g(x)=\bar J(x,\bar\al_\Barr)$.
\end{theorem}
{\bf Proof:}
Recall that $U=\inf_{\bar \A_{\cdot}}\bar J(\cdot,\bar\al)$. Therefore, $J (\cdot,\bar\al_\Barr)\ge U$.
If we show that $g\ge\bar J (\cdot,\bar\al_\Barr)$ then together with Theorem \ref{thm_viscosity_g} it will follow that $U= g\ge\bar J (\cdot,\bar\al_\Barr)$, and so
$U= g=\bar J (\cdot,\bar\al_\Barr)$.
Hence in what follows we shall prove that $g\ge \bar J(\cdot,\bar\al_\Barr)$.
By the definition of the $\Barr$-barrier strategy and the function $g$ it follows that for every $x\in[\Barr,D]$ one has $g(x)=g(\Barr)+r(x-\Barr)$ and $\bar J(x,\bar\al_\Barr)=\bar J(\Barr,\bar\al_\Barr)+r(x-\Barr)$. Therefore, it is sufficient to prove that $g(x)\ge \bar J(x,\bar\al_\Barr)$ only for $x\in[0,\Barr]$.
Apply the definition of $\bar J$ for $\bar\al=\bar\al_\Barr$, that is,
\[
\bar J(x,\bar\al_\Barr)=\sup_{ \psi\in \calQ[x,\bar\al_\Barr]} \bar c(x,\psi,\bar\al_\Barr[\psi]).
\]
Then it is sufficient to show that for every $x\in[0,\beta_0]$ and every $\psi\in \calQ[x,\bar\al_\Barr]$ one has
$$\bar c(x,\psi,\bar\al_\Barr[\psi]) = \int_0^\tau [-\cco(\dot\psi^1)^2(t)-\cct(\dot\psi^2)^2(t) + h(\ph(t))]dt+r\vr(\tau)\leq g(x),$$ 
where $\vr=\bar\al_\Barr[\psi]$ and $\ph(t)=x+yt+\psi^1(t)-\psi^2(t)-\vr(t)$, $0\le t\le\tau$.

Fix such $x$ and $\psi$.
Applying the chain rule to $g$ gives\footnote{Notice that since $x\le \Barr$, and by the definition of the $\Barr$-barrier strategy the dynamics admit no jumps.}
\begin{align}\label{eq085}
&\int_0^\tau [-\cco(\dot\psi^1)^2(t)-\cct(\dot\psi^2)^2(t) +h(\ph(t))]dt + r\vr(\tau)\\\notag
&\quad\leq \int_0^\tau [-\cco(\dot\psi^1)^2(t)-\cct(\dot\psi^2)^2(t) +h(\ph(t))]dt + r\vr(\tau) + g(\ph(\tau))\\\notag
&\quad=g(x) +  \int_0^\tau [-\cco(\dot\psi^1)^2(t)-\cct(\dot\psi^2)^2(t) +h(\ph(t)) + (\dot\psi^1(t)-\dot\psi^2(t)+y)g'(\ph(t))]dt \\\notag
&\quad\quad\qquad\,+ \int_0^\tau[r-g'(\ph(t))]d\vr(t)\\\notag
&\quad\leq g(x) +  \int_0^\tau[r-g'(\Barr)]d\vr(t),
\end{align}
where we write $g'(\Barr)$ for $g_L'(\Barr)$.
The last inequality above holds due to the fact that for every $0\le t\le\tau$ one has $\ph(t)\in[0,\Barr]$ and
$$
\max_{p,q\in\R}\left[-\cco p^2(t)-\cct q^2+h(\ph(t)) + (p-q+y)g'(\ph(t))\right] = -\calL(g'(\ph(t))+h(\ph(t)),
$$ 
and the r.h.s.\ above equals zero by the choice of $g$.

Next, recall that the measure $d\vr$ charges only times $t$ at which $\{\ph(t)=\Barr\}$.
We analyze the last term in \eqref{eq085}
separately for $r/(2\cc) +y \leq 0$ and $r/(2\cc) +y > 0$.
Recall by Section \ref{sec_3c} that in case $r/(2\cc) +y \leq 0$, the function $g$ is differentiable at $\Barr$ and satisfies $g'(\Barr)=r$. Therefore, $\int_0^\tau[r-g'(\Barr)]d\vr(t)=0$ and from \eqref{eq085} it follows that $g(x)\geq \bar J(x,\bar\al_\Barr)$.

In case $r/(2\cc) +y >0 $ the above argument does not hold since $g_L'(\Barr)<g_R'(\Barr)=r$.
Hence we refine the proof.
We will construct a sequence of functions $\{\gd\}_\delta\subset\calC^1([0,D],\R_+)$, that converges uniformly to $g$ as $\delta\rightarrow 0$ with the property that for every
$x\in[0,\Barr]$ and every $\psi\in\calQ[x,\bar\al_\Barr]$ one has
\begin{equation}\label{eq087}
\int_0^\tau [-\cco(\dot\psi^1)^2(t)-\cct(\dot\psi^2)^2(t)+h(\ph(t))]dt + r\vr(\tau)\leq\gd(x),
\end{equation}
by which the result thus follows upon taking $\delta\rightarrow 0$.

\begin{lemma}\label{lem_gd}
There exists a sequence of functions $\{\gd\}_\delta$ that satisfies the following conditions\\
\noi (c1) $\{\gd\}_\delta\subset\calC^1([0,D],\R_+)$.\\
(c2) $\gd'(\Barr)=g_R'(\Barr)=r$.\\
(c3) For every $x\in[0,\Barr]$ one has $\calL(\gd'(x))-h(x)\ge 0$.\\
(c4) $\lim_{\delta\rightarrow 0}\sup_{x\in[0,D]}|\gd(x)-g(x)|=0$.
\end{lemma}
The proof appears below. To complete the proof of the theorem,
note by condition (c1) that one can apply the chain rule to $\gd$, and as in \eqref{eq085},
\begin{align}\notag
&\int_0^\tau [-\cco(\dot\psi^1)^2(t)-\cct(\dot\psi^2)^2(t)+h(\ph(t))]dt + r\vr(\tau)\\\notag
&\quad\leq \int_0^\tau [-\cco(\dot\psi^1)^2(t)-\cct(\dot\psi^2)^2(t)+h(\ph(t))]dt + r\vr(\tau) + \gd(\ph(\tau))\\\notag
&\quad=\gd(x) +  \int_0^\tau [-\cco(\dot\psi^1)^2(t)-\cct(\dot\psi^2)^2(t)+h(\ph(t)) + (\dot\psi^1(t)-\dot\psi^2(t)+y)\gd'(\ph(t))]dt \\\notag
&\quad\quad\qquad\;\;+  \int_0^\tau[r-\gd'(\ph(t))]d\vr(t)\\\notag
&\quad\leq \gd(x) +  \int_0^\tau [-\calL(\gd'(\ph(t)))+h(\ph(t))]dt+  \int_0^\tau[r-\gd'(\Barr)]d\vr(t).
\end{align}
From conditions (c3) and (c2) we get that the second term on the r.h.s.\ is nonpositive and the third one equals zero, and therefore \eqref{eq087} holds.
\hfill $\Box$

\skp\noi{\bf Proof of Lemma \ref{lem_gd}:}
 In this proof we write $g'(\Barr)$ for $g_L'(\Barr)$.
 Fix $0<\delta<\Barr$. Recall that the function $g$ is not differentiable (only) at $\Barr$ and satisfies $g_L'(\Barr)<g_R'(\Barr)=r$.
The idea of the proof is to construct a smooth function $\gd$ that is close to $g$
and such that $\gd=g$ on $[0,\xd]$ where $\xd$ is a little smaller than $\Barr$. The function $g_\delta$ increases on the interval $[\xd,\Barr]$ and satisfies
$\gd'(\Barr)=g_R'(\Barr)=r$.
On $[\Barr,D]$ it is affine with slope $r$.
We start with a linear interpolation. Let $\ld$ be the linear interpolation between the two-dimensional points $(\Barr-\delta, g'(\Barr-\delta))$ and $(\Barr,g_R'(\Barr))$. That is
$\ld(x)=g'(\Barr-\delta)+(x-(\Barr-\delta))\frac{g_R'(\Barr)-g'(\Barr-\delta)}{\delta}$. Define $\xd:=\arg\max\{x\in[\Barr-\delta,\Barr] : \ld(x)= g'(x)\}$. Since $\ld$ and $g'$ are continuous and since $\ld(\Barr-\delta)=g'(\Barr-\delta)$ and $\ld(\Barr)>g_L'(\Barr)$ it follows that $\Barr-\delta\leq\xd<\Barr$.
Define
\[
\gd(x):=\left\{\begin{array}{ll}\notag
               g(x) &\ 0\leq x \leq \xd,
\\
               g(\xd) + \int_{\xd}^x \ld(u)du  & \   \xd <x\leq \Barr,
\\
                \gd(\Barr)+r(x-\Barr) & \  \Barr< x \leq D.
              \end{array}
\right.
\]
One can easily verify that (c1) and (c2) hold. Condition (c3) holds trivially on $[0,\xd]$ since on this interval $\gd=g$ and we chose $g$ to satisfy $\calL(g'(x))-h(x) = 0$ on $[0,\Barr)$. We now turn to the interval $(\xd,\Barr]$. Fix $x\in(\xd,\Barr]$. According to the definition of $\xd$ one has $\gd'(x)=\ld(x) > g'(x)$. Also, $\gd'(x)\leq g_R'(\Barr)=r$. Altogether,
\begin{align}\label{eq090}
g'(x)< g_\delta'(x)\le r.
\end{align}
We claim that $\calL(\gd'(x))-h(x)\ge 0$. To this end, we show that
\begin{align}\label{eq091}
\calL(g'(x))-h(x)= 0\quad\text{and}\quad \calL(r)-h(x)\ge 0,\quad x\in(\xd,\Barr]
\end{align}
and use inequality \eqref{eq090} and the concavity of $\calL$ to deduce that (c3) holds.
We already showed that the first condition of \eqref{eq091} holds. As for the second, notice that since $h$ is increasing, one has $\calL(r)-h(x)\ge \calL(r)-h(\Barr)\ge 0$, where the last inequality follows by the definition of $\Barr$.

To verify (c4) it is sufficient to show that
$\lim_{\delta\rightarrow 0}\sup_{x\in[\xd,\Barr]}|\gd(x)-g(x)|=0$.
Since $\gd'(x)-g'(x)=\ld(x)-g'(x)>0$ for every $x\in[\xd,\Barr]$ it follows that $\gd(x)-g(x)$ increases on $[\xd,\Barr]$, and therefore
we need only to prove that
$\lim_{\delta\rightarrow 0}(\gd(\Barr)-g(\Barr))=0$. This condition indeed holds, due to the fact that $\gd(\Barr) =g(\xd) + \int_{\xd}^\Barr \ld(u)du$, $\lim_{\delta\rightarrow0}\xd=\Barr$, and $\ld$ is bounded.
\hfill $\Box$

As an immediate consequence we obtain an analogous result for the original game.

\begin{theorem}\label{thm_g>J_2}
Let $\al_\Barr$ be the $\Barr$-barrier strategy in the original game. For every $x\in[0,D]$ one has $g(x)= J(x,\al_\Barr):=\sup_{ \psi\in \calP^2, T\in\R_+} c(x,T,\psi,\al_\Barr[\psi])$.
\end{theorem}
{\bf Proof:}
In the proof of Theorem \ref{thm_g>J}, replace the set $\calQ[x,\bar\al_\Barr]$ by $\calP^2$ and the hitting time $\tau$ by any time $T\in\R_+$ to obtain
$$
g(x)= J(x,\al_\Barr)=\sup_{ \psi\in \calP^2, T\in\R_+} c(x,T,\psi,\al_\Barr[\psi]).
$$
\qed

\subsection{Further properties of the game}\label{sec_3d}

This section is concerned with two particular controls for the maximizer, that besides
being interesting by their own right, play an important role in the analysis of the asymptotic stochastic control problem in the forthcoming paper \cite{atar-cohen}.
The first control regards the behavior of the maximizer when the
dynamics are above the free boundary point $\beta_0$, and is used to show
that the property of the $\beta_0$-barrier policy,
by which the dynamics must initially jump from $x>\beta_0$ to $\beta_0$,
is necessary for optimality of the strategy (Subsection \ref{subsec1}).
The second control presented (Subsection \ref{subsec2})
is the optimal response to the $\beta_0$-barrier strategy.
Finally, in the case when $h$ is convex,
the latter control is shown (in Subsection \ref{subsec3})
to correspond to a saddle point for the game.
Specifically, by choosing this control, the maximizer can force a payoff
of at least $U(x)$ regardless of the strategy used by the minimizer.

\subsubsection{Dynamics outside the free boundary}\label{subsec1}
For initial condition $x>\beta_0$,
the optimal strategy that we have identified dictates that the dynamics
make a jump to the point $\beta_0$ at time zero. The initial (i.e.,
at time zero)
cost associated with this jump is $r(x-\beta_0)$.
We will show that any strategy that does not make such an initial jump
is suboptimal. We do so by identifying a control
that, if used by the maximizing player, assures that under any
strategy that does not perform such a jump,
a larger cost is incurred up to the first time $\beta_0$ is reached.

The precise formulation of this requires an additional parameter,
$\del>0$, that can be made arbitrarily small.
Fix $x>\Barr$. Let $\delta>0$ be such that $x>\Barr+\delta$. 
Fix a strategy $\bar\al\in\bar\A_x$.
Consider the
control\footnote{Notice that $\psi_\delta\in\calQ[x,\bar\al]$. }
\begin{equation}
\psi(t)=\psi_\del(t):=\left\{\begin{array}{ll}
               (rt/(2\cco),-rt/(2\cct)) &\  0\le t< T,
\\
                (rT/(2\cct),rT/(2\cct)) &\  t\ge T,
              \end{array}
\right.
\end{equation}
where
$T=T_\delta:= \frac{r(D-\Barr)}{h(\Barr+\delta)-h(\Barr)}$.
Denote also $\ph_{\bar\al}:=\ph[x,\psi,\bar\al[\psi]]$
and let $$\tau=\tauald:= \inf\{t\ge 0 : \ph_{\bar\al}(t)\leq\Barr+\delta \}$$ 
be the first time that the dynamics cross $\Barr+\delta$.

\begin{proposition}\label{prop_V=VB+r(x-B)}
If $\tau>0$ then one has
\begin{align}\notag
\int_0^{\tau}h(\ph_{\bar\al}(t))dt+r\bar\al[\psi](\tau)-\Ir(\tau,\psi)> r(x-(\Barr+\delta)).
\end{align}
Notice that the l.h.s.~of the above is the cost incurred up to time $\tau$, whereas
the r.h.s.\ corresponds to cost of rejecting at time zero.
\end{proposition}

\skp\noi{\bf Proof of Proposition \ref{prop_V=VB+r(x-B)}:}
We analyze separately the cases $0<\tau\le T$ and $\tau>T$. In case $\tau<T$
the cost until time $\tau$ is given by
\begin{align}\notag
&\int_0^{\tau}[-\cco(\dot\psi^1)^2(t) -\cct(\dot\psi^2_\delta)^2(t) +h(\ph(t))]dt + r\vr(\tau) \\\label{eq062}
&\quad\ge
(h(\Barr+\delta)-\frac{r^2}{4\cc})\tau + r\vr(\tau)  \\\label{eq063}
&\quad\ge
r(x-\Barr-\delta)+(\frac{r^2}{4c}+ry+h(\Barr))\tau+(h(\Barr+\delta)-h(\Barr))\tau\\\label{eq064}
&\quad=
r(x-\Barr-\delta) +(h(\Barr+\delta)-h(\Barr))\tau.
\end{align}
Above, inequality \eqref{eq062} follows since on the time interval $[0,\tau)$ one has $\ph>\Barr+\delta$ and since $h$ is increasing. Inequality \eqref{eq063} follows since
$\Barr+\delta\ge\ph(\tau) = x+(\frac{r}{2\cc}+y)\tau - \vr(\tau)$. Equation \eqref{eq064} follows since in case $\Barr<D$ one has $r^2/(4\cc)+ry+h(\Barr)=0$ by \eqref{02}.
Since $h$ is increasing and $\tau>0$, the expression in \eqref{eq064} is strictly greater than $r(x-(\Barr+\delta))$.

In case $\tau>T$ one has $\dot\psi(t)=(0,0)$ on the time interval $[T,\tau]$. Therefore,
the cost incurred until $\tau$ is at least $\int_0^{T}[-\cco(\dot\psi^1)^2(t) -\cct(\dot\psi^2_\delta)^2(t) +h(\ph(t))]dt + r\vr({T})$.
By similar arguments one has that
\begin{align}\notag
&\int_0^T[-\cco(\dot\psi^1)^2(t) -\cct(\dot\psi^2)^2(t) +h(\ph(t))]dt + r\vr(T) \notag
\ge
(h(\Barr+\delta)-\frac{r^2}{4\cc})T + r\vr(T) \\\notag
&\quad\ge
r(x-\Barr-\delta)-r(D-\Barr-\delta)+(\frac{r^2}{4c}+ry+h(\Barr))T+(h(\Barr+\delta)-h(\Barr))T\\\notag
&\quad=
r(x-\Barr-\delta) +r\delta>r(x-\Barr-\delta) .
\end{align}

\noi
\hfill$\Box$

\subsubsection{Optimal response to the barrier strategy}\label{subsec2}

Recall that the $\Barr$-barrier strategy is optimal under our standing assumption
$-y\ge r/(4c)$.
We now identify a control for the maximizer
that is optimal when the minimizer uses the $\Barr$-barrier strategy.
More precisely,
in case $y^2-h(\Barr)/\cc>0$ or $x<\Barr$ (resp., $y^2-h(\Barr)/c= 0$ and $x=\Barr$) we find an explicit expression for the maximizer's optimal (resp., $\delta$-optimal)
control and for the termination time.\footnote{One can show that $y^2-h(\Barr)/\cc\ge 0$ and that equality holds if and only if $-y=r/(2\cc)$ and $-y\le r/\cc+h(D)/r$.}

Throughout, it is assumed that $x\in[0,\Barr]$ and that the maximizer uses
the $\beta_0$-barrier strategy, $\bar\al_{\beta_0}$.
We use a constant $\del$ that in
case $y^2-h(\Barr)/\cc>0$ or $x<\Barr$ is set to zero, and in case $y^2-h(\Barr)= 0$
and $x=\Barr$ is assumed to satisfy $0< \delta\le x$.

In this subsection we assume that $h$ is also Lipschitz continuous on $[0,D]$.
Denote $\tilde\tau=\delta+\int_0^{x+y\delta} 1/\sqrt{y^2-h(\xi)/\cc}\;d\xi$ and let $\om\in \calC_0([0,\tilde\tau],\R)\cap\calC^1([\delta,\tilde\tau],\R)$ be defined
as $0$ on $[0,\del]$ and let it be the unique solution of
\begin{equation}\label{eq500}
\dot\om(t)=-y-\sqrt{y^2-h((x+yt+\om(t))\wedge(x+y\delta))/\cc},\quad t\in[\delta,\tilde\tau].
\end{equation}
Indeed, existence and uniqueness will follow from the Picard-Lindel{\"o}f theorem
once we show that that the r.h.s.~of \eqref{eq500} is Lipschitz with
respect to $\om$. This follows since $h$ is assumed Lipschitz,
and since the expression under the root sign is bounded away from zero.
For the last statement, recall that we assumed that $\delta=0$ in case $y^2-h(\Barr)/\cc>0$ or $x<\Barr$ and $0<\delta<x$ in case $y^2-h(\Barr)/\cc=0$ and $x=\Barr$. Hence, $y^2 - h(x+y\delta)/\cc>0$. Using moreover that $h$ is increasing we get that $y^2-h((x+yt+\om(t))\wedge(x+y\delta))/\cc$ is bounded away from zero.

Set
\begin{align}\label{eq500a}
\psi(t):=\left(\frac{\cc}{\cco}\om(t),\frac{-\cc}{\cct}\om(t)\right),\quad t\in[0,\tilde\tau].
\end{align}

\begin{theorem}\label{thm_501}
If $y^2-h(\Barr)/\cc>0$ or $x<\Barr$ then $\bar c(x, \psi, \bar\al_\Barr[\psi]) = U(x)$.
If $y^2-h(\Barr)/\cc=0$ and $x=\Barr$ then there is a positive constant $C$ such that for every $0<\delta\le x$ one has $|\bar c(x, \psi, \bar\al_\Barr[\psi])-U(x)|\le C\delta$.
Moreover, for every $\delta\ge 0 $ the termination time of the game is $\tau=\tilde\tau$.
\end{theorem}

\skp\noi{\bf Proof:}
Notice that by the optimality of $\bar\al_\Barr$  it follows that for every $\delta\ge 0$ one has $U(x)=\sup_{\tilde\psi\in\calQ[x,\bar\al_\Barr]}\bar c(x, \tilde\psi, \bar\al_\Barr[\tilde\psi])\ge \bar c(x, \psi, \bar\al_\Barr[\psi]) $. Therefore, in order to prove the theorem it is sufficient to prove that $\bar c(x,\psi, \bar\al_\Barr[\psi])\ge U(x)-C\delta$ where $\delta=0$ in case $y^2-h(\Barr)/\cc>0$ or $x<\Barr$ and $0<\delta<x$ in case $y^2-h(\Barr)/\cc=0$ and  $x=\Barr$.
To this end, it is sufficient to show that $\psi\in\calQ^2[x,\bar\al_\Barr]$ and that
the cost until the first hitting time, $\tau$, is bounded below by $g$. That is
\begin{align}\label{eq501}
\bar c(x,\psi,\bar\al[\psi]) &=\int_0^{\tau} [-\cco(\dot\psi^{1})^2(t)-\cct(\dot\psi^{2})^2(t)+h(\ph(t))]dt + r\vr(\tau)
\ge g(x),
\end{align}
where $\vr=\bar\al_\Barr[\psi]$ and
$\ph(t)=x+yt+\psi^{1}(t)-\psi^{2}(t)-\vr(t)$, $0\leq t\leq\tau$.
We do that by showing that $\psi$ satisfies the following three conditions.\\
\noi (a) For every $0\leq t\leq\tau$ one has $\ph(t)\in[0,\Barr]$.\\
(b) $\psi\in\calQ^2[x,\bar\al_\Barr]$.\\
(c) For every $\delta\leq t\leq\tau$ one has
\begin{align}\notag
-\cco(\dot\psi^{1})^2(t)-\cct(\dot\psi^{2})^2(t)+h(\ph(t)) + (\dot\psi^{1}(t)-\dot\psi^{2}(t)+y)g'(\ph(t))= 0.
\end{align}

Before proving (a)--(c) let us show that they are sufficient for \eqref{eq501}.
Condition (a) is needed for $g'(\ph(\cdot))$ to be well-defined, where one uses
the convention $g'(\Barr)=g_L'(\Barr)$. From condition (b) it follows that
$g(\ph(\tau))=0$.  As a result, we can apply the chain rule to $g$,
as in \eqref{eq052}. Thus
\begin{align}\notag
&\int_0^{\tau} [-\cco(\dot\psi^{1})^2(t)-\cct(\dot\psi^{2})^2(t)+h(\ph(t))]dt + r\vr(\tau) + g(\ph(\tau))\\\notag
&\quad=g(x) +  \int_0^\delta [-\cco(\dot\psi^{1})^2(t)-\cct(\dot\psi^{2})^2(t)+h(\ph(t)) + (\dot\psi^{1}(t) -\dot\psi^{2}(t)+y)g'(\ph(t))]dt \\\notag
&\quad\quad+ \int_\delta^{\tau} [-\cco(\dot\psi^{1})^2(t)-\cct(\dot\psi^{2})^2(t)+h(\ph(t)) + (\dot\psi^{1}(t) -\dot\psi^{2}(t)+y)g'(\ph(t))]dt\\\notag
&\quad\quad + \int_0^{\tau}[r-g'(\ph(t))]d(\vr)(t).
\end{align}
From condition (c) together with the fact that $ r\ge g'$ it follows that
\begin{align}\label{eq505}
&\int_0^{\tau} [-\cco(\dot\psi^{1})^2(t)-\cct(\dot\psi^{2})^2(t)+h(\ph(t))]dt + r\vr(\tau)\\\notag
&\hspace{.5em}
\ge g(x) + \int_0^\delta [-\cco(\dot\psi^{1})^2(t)-\cct(\dot\psi^{2})^2(t)+h(\ph(t)) + (\dot\psi^{1}(t) -\dot\psi^{2}(t)+y)g'(\ph(t))]dt.
\end{align}
Recall that for every $t\in[0,\delta]$ one has $\dot\om(t)=0$ and therefore $\psi(\cdot)=(0,0)$ on that time interval. Moreover, the functions $h$ and $g$ are both bounded.
Therefore, the integral on the r.h.s.~of \eqref{eq505} is bounded below by $-C\delta$, where $C$ is a positive constant independent of $\delta$.
Hence, we proved that $ c(x, \psi, \bar\al_\Barr[\psi])\ge U(x)-C\delta$.

We now show that $\psi$ satisfies conditions (a)--(c). Hereafter, we assume that $t\ge\delta$. Since $y+\dot\om(t)\leq 0$ it follows that $x+yt+\om(t)$ is nonincreasing  in $t$ and therefore $x+yt+\om(t)\leq x+y\delta\le\Barr$, since moreover, $\om = \psi^{1}-\psi^{2}$ it follows that $\ph(t)=x+yt+\om(t)$ and condition (a) holds. Since $h$ is increasing we get that
\begin{align}\notag
y+\dot\om(t)\leq -\sqrt{y^2-h(x+y\delta)/\cc}<0.
\end{align}
That is, $\ph(\cdot)=x+y\cdot+\om(\cdot)$ is (strictly) decreasing with a slope that is bounded above by a negative number and therefore $\ph$ reaches zero at a finite time and (b) holds.
Condition (c) follows by substitution of $\psi$ and $g$.

We now turn to the termination time of the game. Recall that for every $0\le t\le\tau$ one has $\ph(t)\in[0,\Barr]$. Also, recall that $\ph(\cdot) = x+y\cdot+\om(\cdot)$ and therefore,
$$\tau = \inf\{t\ge 0 : \ph(t)=0 \} = \inf\{t\ge 0 : x+yt+\om(t)=0 \}.$$
As mentioned earlier, the map $\ph: [0,\tau]\to [0,x]$ is strictly decreasing.
Therefore, $\ph$ has an inverse function $\ph^{-1}:[0,x]\to [0,\tau]$. 
Its derivative satisfies
$$(\ph^{-1})'(\xi)=\frac{1}{\ph'(\ph^{-1}(\xi))},$$
and more explicitly
\begin{equation}\label{eq508}
(\ph^{-1})'(\xi):=\left\{\begin{array}{ll}
               -1/\sqrt{y^2-h(\xi)/\cc} \quad  &\quad\quad\; 0\leq \xi \leq x+y\delta,
\\
                   1/y \quad &x+y\delta\leq \xi \leq x.
              \end{array}
\right.
\end{equation}
Since $\ph^{-1}(0)=\tau$ and $\ph^{-1}(x)=0$ it follows that
$$\tau = \ph^{-1}(0)-\ph^{-1}(x) = -\int_0^x (\ph^{-1})'(\xi)d\xi =\delta+ \int_0^{x+y\delta} 1/\sqrt{y^2-h(\xi)/\cc}\;d\xi .$$

\noi
\hfill $\Box$

\subsubsection{The case of convex $h$: a saddle point property}\label{subsec3}

We now present a saddle point property of the game in case when $h$ is a convex
function and twice continuously
differentiable, namely that the maximizer can force a payoff of at least $U(x)$ by
choosing the control $\psi$ defined in the previous subsection.
The result is only concerned with the case $x<\beta_0$. Note that by the convention
regarding the constant $\del$, we take here $\del=0$.
\begin{proposition}\label{prop501}
Fix $x\in[0,\Barr)$ and assume $h$ is a convex function.
Then for every $\vr\in\calE$ one has $\bar c(x,\psi, \vr)\ge U(x)$.
\end{proposition}

\noi{\bf Proof:}
We prove the proposition in two steps. First, we show that if the minimizer uses a control $\vr$ that has $p$ rejections in total then it is optimal for it to
reject the amount $p$ at time zero. Then we argue that, for any amount $p$ rejected at time zero, the cost is bounded below by $U(x)$.

\skp\noi{\bf Step 1.}
Let $\vr\in\calE$. Denote by $p:=\vr(\tau)$, where $\tau=\tau(x,\psi,\vr)$. Define the modified path $\vr^p$ as $\vr^p(t)=p$ for all $t\ge0$. We will show that
\begin{align}\label{eq510}
\bar c(x,\psi,\vr) \ge \bar c(x,\psi,\vr^p).
\end{align}
To this end, recall the definition of $\psi$ from \eqref{eq500a}, which yields that $x+y\cdot+\psi^{1}(\cdot)-\psi^{2}(\cdot)$ is strictly decreasing and therefore
$\tau(x,\psi,\vr)=\tau(x,\psi,\vr^p)=T$,
where
\begin{align}\label{110}
T:=\inf\{t\ge 0 : x+yt+\psi^{1}(t)-\psi^{2}(t)=p\}.
\end{align}
Therefore
\begin{align}\notag
\bar c(x,\psi,\vr)
= \int_0^{T}[-\cco(\psi^{1})^2(t)-\cco(\psi^{2})^2(t) + h(x+yt+\psi^{1}(t)-\psi^{2}(t)-\vr(t))]dt + rp,
\end{align}
and $\bar c(x,\psi,\vr^p)$ is given by a similar expression with $\vr(t)$
replaced by $p$.
Since $h$ is increasing, \eqref{eq510} follows.

\skp\noi{\bf Step 2.}
Denote $\ph(t) = x+yt+\psi^{1}(t)-\psi^{2}(t)$. Note that this is the dynamics
under the control $\psi$ and the $\Barr$-barrier strategy, because $\vr=0$
under that strategy for $x<\beta_0$. Consider any $p\in[0,x]$ and recall
the notation $\vr^p$ from the previous paragraph. Notice that
$\ph[x,\psi,\vr^p]= \ph - p$.
Denote $F(p)= \bar c(x,\psi,\vr^p)$. The result will follow
once we show that $F(p)\ge F(0)=g(x)$. To this end, write $T(p)$ for
the quantity $T$ of \eqref{110}, to emphasize the $p$-dependence. Then
\begin{align}\notag
F(p) &= \int_0^{T(p)}[-\cco(\psi^{1})^2(t)-\cco(\psi^{2})^2(t) + h(\ph(t)-p)]dt + rp\\\notag
       &=\int_0^{T(p)}[-\cco(\psi^{1})^2(t)-\cco(\psi^{2})^2(t) + h(\ph(t))]dt
       \\ \notag
       &\qquad -\int_0^{T(p)}[h(\ph(t))-h(\ph(t)-p)]dt + rp\\\notag
       &= g(x)-g(p) -\int_0^{T(p)}[h(\ph(t))-h(\ph(t)-p)]dt + rp.
\end{align}
The last equality follows since $$g(x) = \int_0^{T(p)}[-\cco(\psi^{1})^2(t)-\cco(\psi^{2})^2(t) + h(x+yt+\psi^{1}(t)-\psi^{2}(t))]dt +g(p),$$
which in turn follows since both sides of the last equation represents the cost when both players play optimally where we use a DPP argument. Indeed, the integral is the cost under the dynamics hits $p$ and $g(p)$ is the cost when the initial state is $p$.

Note that $F(0)=g(x)$.
Thus to prove that $F(p)\ge g(x)$ for every $p\in[0,x]$ it suffices
to show that $F'\ge 0$ on $[0,x]$.
To this end, first notice that $T=\ph^{-1}$ and that $T'(\cdot) = -1/\sqrt{y^2-h(\cdot)/\cc} $, see \eqref{eq508}.
By direct calculations it follows that
\begin{align}\notag
F'(p)=-g'(p)-T'(p)h(p) - \int_0^{T(p)}h'(\ph(t)-p)dt +r.
\end{align}
By substituting $q=\ph(t)$ on the integral on the r.h.s.~of the we get
\begin{align}\notag
\int_0^{T(p)}h'(\ph(t)-p)dt &= \int_x^p h'(q-p)/(\ph)'(\ph^{-1}(q))\;dq = \int_x^p h'(q-p)T'(q)\;dq \\\notag
&= \int_p^x h'(q-p)/\sqrt{y^2-h(q)/\cc}\;dq.
\end{align}
Using moreover the convexity of $h$ it follows that
\begin{align}\notag
F'(p)\ge -g'(p)-T'(p)h(p) - \int_p^x h'(q)/\sqrt{y^2-h(q)/\cc}\;dq +r =-T'(p)h(p) -g'(x)+r \ge 0.
\end{align}
The equality above follows since $h'(q)/\sqrt{y^2-h(q)/\cc}=g''(q)$. The last inequality follows since $r\ge g'(x)$ and since $T'(p)$ is negative.
\noi
\hfill $\Box$

\appendix

\section{Appendix}\beginsec
\manualnames{A}

\subsection*{State space collapse for the multidimensional DG}
The scaled processes $(\tilde A^n,\tilde S^n)$ are assumed to satisfy a {\it moderate deviation
principle}. To express this assumption, let
$\Jr_k, k=1,2,$ be functions on $\calD([0,T],\R^I)$ defined as follows. For
$\upsi^k=(\psi^k_1,\ldots,\psi^k_I)\in \calD([0,T],\R^I)$,
\begin{align}\label{100}
\Jr_k(T,\upsi^k):=\left\{\begin{array}{ll}
               \sum_{i=1}^I c_{i,k}
               \int_0^T(\dot\psi^k_i)^2(s)ds &\ \mbox{if all}\ \psi^k_i\in\AC_0([0,T],\R),
\\
                \infty & \ \mbox{otherwise},
              \end{array}
\right.
\end{align}
where
\begin{align}\notag
c_{i,1}= \frac{1}{2\la_i\sig^2_{i,IA}}\quad\text{and}\quad c_{i,2}=\frac{1}{2\mu_i\sig^2_{i,ST}},\quad i\in\calI.
\end{align}
Let $\Jr(T,\upsi)=\Jr_1(T,\upsi^1)+\Jr_2(T,\upsi^2)$
for $\upsi=(\upsi^1,\upsi^2)\in \calD([0,T],\R^{2I})$.

We show that the multidimensional game can be reduced into a one dimensional game in the sense that they both share the same value. Towards this, we will need to define the components of the one dimensional cost function: the holding cost, rejection cost, and the rate function.
For $w\in\R_+$, denote
\begin{equation}\label{eq1033}
 h(w):=\inf\{\uh\cdot \uxi : \uxi\in \prod_{i=1}^I [0,D_i], \utheta\cdot \uxi=w\}.
\end{equation}
By the convexity of the set $\prod_{i=1}^I [0,D_i]$ it follows that $ h$ is convex. Moreover, $ h(w)\ge0$ for $w\ge 0$ and equality holds if and only if $w=0$. Therefore, $ h$ is
strictly increasing.
Let
\begin{equation}\label{eq1034}
 r:=\min\{\ur\cdot \uq: \uq\in\R_+^I, \utheta\cdot \uq=1\}.
\end{equation}
As was shown in \cite{ata-shi},
\begin{equation}\notag
 r=r_{i^*}\mu_{i^*}:=\min\{r_i\mu_i:i\in\I\}.
\end{equation}
For every $T\in\R_+$ and $\psi=(\psi^1,\psi^2)\in\calP^2$ set $\Ir(T,\psi)=\Ir_1(T,\psi^1)+\Ir_2(T,\psi^2)$, where
\begin{equation}\label{eq1036}
\Ir_k(T,\psi^k) := \inf \{\Jr_k(T,\upsi^k):\upsi^k\in\calP^I, \utheta\cdot\upsi^k=\psi^k\},\quad k=1,2.
\end{equation}
One can verify that
\begin{equation}\notag
\Ir_k(t,\psi^k)=\left\{\begin{array}{ll}
               c_k
               \int_0^T(\dot\psi^k)^2(s)ds &\ \mbox{if}\ \psi^k\in\AC_0([0,T],\R),
\\
                \infty & \ \mbox{otherwise},
              \end{array}
\right.
\end{equation}
where
\begin{equation}\notag
c_1:=\left(\sum_{i=1}^{I}\frac{2\rho_i\sig^2_{i,IA}}{\mu_i}\right)^{-1},\quad\text{and}\quad c_2:=\left(\sum_{i=1}^{I}\frac{2\sig^2_{i,ST}}{\mu_i}\right)^{-1}.
\end{equation}
Moreover, for every $\psi=(\psi^1,\psi^2)\in\calP^2$ define the multidimensional path $\upsi_\psi=(\upsi^1_{\psi^1},\upsi^2_{\psi^2})\in\calP^{2I}$ by
\begin{equation}\notag
\upsi^1_{\psi^1} := \left(\frac{2\rho_1\sig^2_{1,IA}}{\mu_1},\ldots,\frac{2\rho_n\sig^2_{n,IA}}{\mu_n}\right)c_1\psi^1
\end{equation}
and
\begin{equation}\notag
\upsi^2_{\psi^2} := \left(\frac{2\sig^2_{1,ST}}{\mu_1},\ldots,\frac{2\sig^2_{n,ST}}{\mu_n}\right)c_2\psi^2.
\end{equation}
Simple calculation yields that
\begin{equation}\notag
\utheta\cdot\upsi^k_{\psi^k}=\psi^k,\quad k=1,2,
\end{equation}
and that
\begin{equation}\label{eq1042}
\Jr(T,\upsi_\psi)=\Ir(T,\psi).
\end{equation}
Given $\ux\in \prod_{i=1}^I [0,D_i]$,
let $x=\utheta\cdot\ux$. The one-dimensional value function is given by
\begin{equation}\notag
V_{OD}(x)=\inf_{\alpha\in \A_x}\sup_{\psi\in \calP^2, T\in\R_+}\left(\int_0^T h(\varphi(t))dt+r\al_1[\psi](T)-\Jr(T,\psi)\right).
\end{equation}

\begin{theorem}\label{thm1000}
For every $\ux\in \prod_{i=1}^I [0,D_i]$ one has $V_{MD}(\ux)=V_{OD}(x)$.
\end{theorem}
\skp\noi{\bf Proof of Theorem \ref{thm1000}}
Fix $\ux$.
We prove the theorem by showing
\begin{align}\label{eq1044}
V_{OD}( x)\le V_{MD}(\ux)
\end{align}
and
\begin{align}\label{eq1045}
V_{MD}(\ux)\ge V_{OD}(x).
\end{align}
We start by showing \eqref{eq1044}.
Let $\ual\in\A_{\ux}$. Define $\al\in\A_x$ by
\begin{align}\notag
\al_k[(\psi^1,\psi^2)]:=\utheta\cdot\ual_k[(\upsi^1_{\psi^1},\upsi^2_{\psi^2})],\quad (\psi^1,\psi^2)\in\calP^2,\quad k=1,2.
\end{align}
The proof that $\ual$ is admissible is straightforward and therefore omitted.
Since $(\psi^1,\psi^2) = (\utheta\cdot\upsi^1_{\psi^1},\utheta\cdot\upsi^2_{\psi^2}) $ we get that
\begin{align}\label{eq1047}
c(x,T,\psi, \al[\psi])=c(x,T,(\utheta\cdot\upsi^1_{\psi^1},\utheta\cdot\upsi^2_{\psi^2}), \al[(\utheta\cdot\upsi^1_{\psi^1},\utheta\cdot\upsi^2_{\psi^2})])
\end{align}
Fix $\psi\in\calP^2$ and let
\begin{align}\notag
\uph(t)=\ux+\uy t+\upsi^1_{\psi^1}(t)-\upsi^2_{\psi^2}(t)+\ual_1[\upsi_\psi](t)
-\ual_2[\upsi_\psi](t),\qquad t\ge 0
\end{align}
be the dynamics associated with $\ux,\upsi_\psi$, and $\ual[\upsi_\psi]$. By the definitions of $ h$ and $ r$ (see \eqref{eq1033}, \eqref{eq1034}) and by \eqref{eq1042} it follows that
\begin{align}\notag
\uh\cdot\uph(t)\ge h(\utheta\cdot\uph(t)),\quad t\ge0,
\end{align}
\begin{align}\notag
\ur\cdot \ual_2[\upsi_\psi](t)\ge  r \utheta\cdot \ual_2[\upsi_\psi](t),\quad t\ge0,
\end{align}
and
\begin{align}\notag
\Ir(T,\psi)=\Jr(T,\upsi_\psi),\quad T\ge0.
\end{align}
Therefore,
\begin{align}\notag
&c(x,T,(\utheta\cdot\upsi^1_{\psi^1},\utheta\cdot\upsi^2_{\psi^2}), \al[(\utheta\cdot\upsi^1_{\psi^1},\utheta\cdot\upsi^2_{\psi^2})])\\\notag
&\quad=
\int_0^T h(\utheta\cdot\uph(t))dt + r\utheta\cdot \ual_2[\upsi_\psi](T)-\Ir(T,\upsi_\psi)\\\notag
&\quad\le
\int_0^T \uh\cdot \uph(t)dt + \ur\cdot \ual_2[\upsi_\psi](T)-\Jr(T,\upsi)\\\notag
&\quad=
\uc( \ux,T,\upsi_\psi, \ual[\upsi_\psi]).
\end{align}
Together with \eqref{eq1047} we get that
\begin{align}\notag
c(x,T,(\psi), \al[\psi])\le \uc( \ux,T,\upsi_\psi, \ual[\upsi_\psi])
\end{align}
and therefore,
\begin{align}\notag
\sup_{\psi\in \calP^2, T\in\R_+}c(x,T,(\psi), \al[\psi])
\le
\sup_{\upsi\in \calP^{2I}, T\in\R_+}c( \ux,T,\upsi, \ual[\upsi]).
\end{align}
This shows \eqref{eq1044}.

We next show \eqref{eq1045}. Let $\al\in\A_x$. Toward constructing
$\ual\in\A_{\ux}$ let $\gamma:[0,D]\to \prod_{i=1}^I [0,D_i]$ be Borel measurable, satisfying
\begin{equation}\notag
\gamma(w)\in\argmin_{\uxi}\{\uh\cdot\uxi:\uxi\prod_{i=1}^I [0,D_i],\utheta\cdot\uxi=w\},\qquad w\in[0,D].
\end{equation}
The existence follows from a measurable selection argument
(such as Corollary 10.3 in the appendix of \cite{ethkur}).
Set
\begin{equation}\notag
\ph(t)=\utheta\cdot \ux+\utheta\cdot \uy t +\utheta\cdot\upsi^1(t)-\utheta\cdot\upsi^2(t)+\ual_1[(\upsi^1,\upsi^2)](t)-\ual_2[(\upsi^1,\upsi^2)](t),\quad t\ge 0,\quad (\upsi^1,\upsi^2)\in\calP^{2I}
\end{equation}
and
\begin{equation}\notag
\uph(t)=\gamma(\ph(t)),\quad t\ge 0.
\end{equation}
Define
\begin{equation}\notag
\ual_2[(\upsi^1,\upsi^2)](t):=\al_2[(\utheta\cdot\upsi^1,\utheta\cdot\upsi^2)](t)\mu_{i^*}e_{i^*},\quad t\ge 0,\quad (\upsi^1,\upsi^2)\in\calP^{2I},
\end{equation}
and
\begin{equation}\notag
\ual_1[(\upsi^1,\upsi^2)](t):=\uph(t)-\ux-\uy t -\upsi^1(t)+\upsi^2(t)+\ual_2[(\upsi^1,\upsi^2)](t),\quad t\ge 0,\quad (\upsi^1,\upsi^2)\in\calP^{2I}.
\end{equation}
The strategy $\ual$ is admissible. Indeed, causality follows since $\al$ satisfies this property. Monotonicity of each coordinate of $\ual_2[(\upsi^1,\upsi^2)](\cdot)$ follows by the monotonicity of $\al_2[(\utheta\cdot\upsi^1,\utheta\cdot\upsi^2)](\cdot)$, and monotonicity of $\utheta\cdot\ual_1[(\upsi^1,\upsi^2)](\cdot)$ follows since the latter equals to $\al_1[(\utheta\cdot\upsi^1,\utheta\cdot\upsi^2)](\cdot)$ which is monotone by
admissibility of $\al_1$. Finally, for every $w\in[0,D]$ one has $\gamma(w)\in \prod_{i=1}^I [0,D_i]$ and therefore $\uph$ takes values in $\prod_{i=1}^I [0,D_i]$.

Next, note that
\begin{equation}\notag
\uh\cdot \uph(t)= h( \ph(t)),\quad t\ge 0,
\end{equation}
\begin{equation}\notag
\ur\cdot \ual_2[(\upsi^1,\upsi^2)](t)=r_{i^*}\mu_{i^*}\al_2[(\utheta\cdot\upsi^1,\utheta\cdot\upsi^2)](t)=  r\al_2[(\utheta\cdot\upsi^1,\utheta\cdot\upsi^2)](t),\quad t\ge 0,\quad (\upsi^1,\upsi^2)\in\calP^{2I}
\end{equation}
and
\begin{equation}\notag
\Jr(T,(\upsi^1,\upsi^2))\ge\Ir(T,(\utheta\cdot\upsi^1,\utheta\cdot\upsi^2)),\quad T\ge 0,\quad (\upsi^1,\upsi^2)\in\calP^{2I}.
\end{equation}
The last inequality follows by the definition of $\Ir$, see \eqref{eq1036}.
Therefore for every $\upsi=(\upsi^1,\upsi^2)\in\calP^{2I}$ one has
\begin{align}\notag
&\uc(\ux,T,\upsi, \ual[\upsi])\\\notag
&\quad
=
\int_0^T \uh\cdot \uph(t)dt + \ur\cdot \ual_2[(\upsi^1,\upsi^2)](T)-\Jr(T,(\upsi^1,\upsi^2))\\\notag
&\quad
\le
\int_0^T h(\ph(t))dt + r\al_2[(\utheta\cdot\upsi^1,\utheta\cdot\upsi^2)](T)-\Ir(T,(\utheta\cdot\upsi^1,\utheta\cdot\upsi^2))\\\notag
&\quad
=
c( x,T,(\utheta\cdot\upsi^1,\utheta\cdot\upsi^2), \al[(\utheta\cdot\upsi^1,\utheta\cdot\upsi^2)])
\end{align}
Hence,
\begin{align}\notag
&\sup_{\upsi\in \calP^{2I}, T\in\R_+}\uc(\ux,T,\upsi, \ual[\upsi])\le
\sup_{\psi\in \calP^2, T\in\R_+}c( x,T,\psi, \al[\psi])
\end{align}
and \eqref{eq1045} follows.

\footnotesize

\bibliographystyle{is-abbrv}


\end{document}